\documentclass[reqno]{amsart}
\usepackage{amsfonts,amssymb,amsmath,mathtools,mathrsfs}
\usepackage{lmodern}
\usepackage[english]{babel}
\usepackage[utf8]{inputenc}
\usepackage[a4paper,margin=2.1cm]{geometry}
\usepackage{xcolor}
\usepackage{cite}
\definecolor{RefBlue}{RGB}{0,0,180}
\usepackage[colorlinks=true,linkcolor=RefBlue,citecolor=RefBlue,
 urlcolor=RefBlue,bookmarksnumbered=true,bookmarksopen=true,
 linktoc=all,breaklinks=true,pdfborder={0 0 0},hypertexnames=false]{hyperref}
\newtheorem{maintheorem}{Theorem}

\newtheorem{theorem}{Theorem}[section]
\newtheorem{proposition}[theorem]{Proposition}
\newtheorem{lemma}[theorem]{Lemma}
\newtheorem{corollary}[theorem]{Corollary}
\theoremstyle{definition}

\theoremstyle{remark}
\newtheorem{remark}[theorem]{Remark}
\numberwithin{equation}{section}
\DeclareMathOperator{\Div}{div}
\DeclareMathOperator{\diag}{diag}
\DeclareMathOperator{\dist}{dist}
\newcommand{\eps}{\varepsilon}
\newcommand{\Om}{\Omega}
\newcommand{\Omeps}{\Omega_\eps}
\newcommand{\Ystar}{Y^*}
\newcommand{\Qeps}{Q_\eps}
\newcommand{\R}{\mathbb R}
\newcommand{\dd}{\,d}
\newcommand{\weakto}{\rightharpoonup}
\begin{document}
\renewcommand{\abstractname}{\normalfont\bfseries Abstract}
\title[Degenerate parabolic Robin homogenization]{Homogenization of degenerate parabolic equations with nonlinear Robin reactions in perforated domains}
\author[Hermann Douanla]{Hermann Douanla$^{*}$}
\author{Rodrigue Kenne Bogning$^{\dagger}$}
\author{Erick Tetsadjio$^{\ddagger}$}
\thanks{$^{*}$Corresponding author. Department of Mathematics,
University of Yaounde I, P.O. Box 812, Yaounde, Cameroon.
\textit{E-mail:} \href{mailto:hermann.douanla@facsciences-uy1.cm}{\texttt{hermann.douanla@facsciences-uy1.cm}}.}
\thanks{$^{\dagger}$Department of Fundamental and Transversal Sciences,
National Advanced School of Public Works, P.O. Box 510, Yaounde, Cameroon.
\textit{E-mail:} \href{mailto:rodriguekenne@ymail.com}{\texttt{rodriguekenne@ymail.com}}.}
\thanks{$^{\ddagger}$Department of Mathematics, University of Yaounde I,
P.O. Box 812, Yaounde, Cameroon.
\textit{E-mail:} \href{mailto:tetsadjio@gmail.com}{\texttt{tetsadjio@gmail.com}}.}
\date{}
\begin{abstract}
We study periodic homogenization of a parabolic equation with linear,
degenerate bulk terms in a perforated domain, with a monotone nonlinear Robin reaction of order
$\varepsilon$ on the pore walls. The storage coefficient may vanish on sets
of positive measure, and the diagonal diffusion weights are not bounded
below but satisfy an inverse-integrability condition. Under positivity of the mean storage, a quantitative cell-detection estimate
yields strong compactness of the states and of the scaled traces. The
effective equation has positive storage, a positive-definite diffusion
matrix given by a weighted cell problem, and a nonlinear volume reaction
obtained by averaging over the pore surface. We also prove a strong
weighted-gradient corrector and uniform-in-time convergence of the
dissipated and stored energies and of the states in the storage-weighted
$L^2$ norm.

\vspace{1.5ex}
\noindent\textbf{2020 Mathematics Subject Classification.} 35B27, 35K20, 35K57, 35K65, 35R05, 76S05.

\noindent\textbf{Key words and phrases.} Periodic homogenization, vanishing storage, degenerate diffusion, nonlinear Robin reaction, perforated domain.
\end{abstract}
\maketitle

\section{Introduction}\label{sec:introduction}

Diffusion through periodically perforated porous materials, coupled with
reactions on the pore walls, arises in heterogeneous catalysis, in
chromatography and in the transport of solutes in soils, where adsorption at
the solid surface is typically described by a nonlinear, monotone isotherm.
In many such materials the medium is far from uniform at the pore scale. Parts
of the solid phase may have negligible storage capacity, so that the
equation is parabolic in some regions and elliptic in others
\cite{AltLuckhaus1983}; and the mobility may deteriorate in some directions,
for instance close to the pore walls.
The present paper studies the effective evolution and the nonlinear reaction
in this coupled degenerate regime.

\subsection{The model}\label{sec:model}

Let $N\geq2$ be an integer, $T>0$, and let $Y=(0,1)^N$ be the reference cell. The porous
medium is obtained from a bounded domain $\Om$ by removing the
$\eps$-scaled periodic copies of a reference hole contained in $Y$, whose
boundaries are the pore walls. We assume the following hypothesis on the
geometry.

\noindent\textbf{(H$_1$).}\enspace
The domain $\Om\subset\R^N$ is bounded and Lipschitz;
$\varnothing\ne \mathcal O\Subset Y$ is open with Lipschitz boundary
$S=\partial \mathcal O$; and $\Ystar=Y\setminus\overline{\mathcal O}$ is connected.

\medskip
For $0<\eps<1$, we perforate only the cells contained in $\Om$ and set
\begin{align*}
 I_\eps&=\{k\in\mathbb Z^N:\eps(k+\overline Y)\subset\Om\},
 &\mathcal O_\eps&=\bigcup_{k\in I_\eps}\eps(k+\mathcal O),\\
 \Omeps&=\Om\setminus\overline{\mathcal O_\eps},
 &S_\eps&=\bigcup_{k\in I_\eps}\eps(k+S),
\end{align*}
together with $Q=\Om\times(0,T)$, $\Qeps=\Omeps\times(0,T)$ and the
unperforated boundary strip
$\Sigma_\eps=\Om\setminus\bigcup_{k\in I_\eps}\eps(k+Y)$, whose measure is
$O(\eps)$ by Lemma~\ref{lem:periodic-averaging}. Since
$\eps^N\#I_\eps\leq|\Om|$,
\begin{equation}
 \eps|S_\eps|=\eps^N\#I_\eps\,|S|\leq |\Om|\,|S|.
 \label{eq:surface-mass}
\end{equation}
The normal $\nu_\eps$ on $S_\eps$ is exterior to $\Omeps$ and points
into the holes.

The coefficients $b$, $A=\diag(w_1,\ldots,w_N)$ and $h$ are
periodic in the microscopic variable and satisfy the conditions in
Section~\ref{sec:hypotheses}. With $b^\eps(x)=b(x/\eps)$,
$w_i^\eps(x)=w_i(x/\eps)$, $\omega^\eps=(w_1^\eps,\ldots,w_N^\eps)$ and
$A^\eps=\diag(w_1^\eps,\ldots,w_N^\eps)$, the microscopic problem is
\begin{equation}
 \left\{
 \begin{aligned}
  \partial_t(b^\eps u_\eps)-\Div(A^\eps Du_\eps)&=f
       &&\text{in }\Qeps,\\
  A^\eps Du_\eps\cdot\nu_\eps
       +\eps h(x/\eps,u_\eps)&=0
       &&\text{on }S_\eps\times(0,T),\\
  u_\eps&=0&&\text{on }\partial\Om\times(0,T),\\
  b^\eps u_\eps(0)&=b^\eps u^{\rm in}&&\text{in }\Omeps.
 \end{aligned}\right.
 \label{eq:micro}
\end{equation}
The factor $\eps$ compensates the growth of the total area of the holes,
and the initial condition is understood as an initial mass in the dual of the
energy space (Section~\ref{sec:micro}). Throughout, $D$ denotes the spatial
gradient and $D_i=\partial/\partial x_i$; for functions of $(x,t,y)$ we write
$D_x,D_y$ and $D_{x_i},D_{y_i}$.

Three difficulties interact. First, $b$ may vanish on sets of positive
measure, so that the time derivative controls only the mass $b^\eps u_\eps$
and the associated mass operator is not injective. The classical
Aubin--Lions argument based on a time-derivative estimate for the state
therefore does not apply directly. Second, the weights
$w_i$ are positive almost everywhere but need not be bounded below; under an
inverse-integrability condition the uniform estimates control the states only
in $W^{1,q}$ with $q<2$, so that neither the $H^1$ regularity of the limit nor
the two-scale structure of the gradients is automatic. Third, the nonlinear
reaction acts on an oscillating surface of total area of order $\eps^{-1}$,
and its limit requires strong convergence of the scaled traces, which has to
be derived from weighted bulk estimates only. Beyond identifying the limit,
we seek estimates uniform in the observation time, using only the time
regularity of the mass.

\subsection{Related work}

\emph{Boundary conditions on periodic holes.}
Nonhomogeneous Neumann and Robin conditions on holes of size $\eps$ go back to
\cite{CioranescuDonato1988,Conca1985}. The periodic unfolding method
\cite{CioranescuDamlamianGriso2018}, adapted to perforated domains in
\cite{CioranescuDamlamianDonatoGrisoZaki2012}, was applied to linear and
nonlinear Robin problems in
\cite{CioranescuDonatoZaki2006,CioranescuDonatoZaki2007}, and to
quasilinear problems with nonlinear Robin conditions and $L^1$ data in
\cite{CabarrubiasDonato2012,DonatoGuibeOropeza2018}; corrector results in
this stationary setting were obtained through truncations of renormalized
solutions \cite{DonatoGuibeOropeza2023}. Correctors for the heat equation
in perforated domains go back to \cite{DonatoNabil2001}. Surface terms can
also be handled by two-scale convergence on periodic surfaces
\cite{AllaireDamlamianHornung1996,NeussRadu1996}.

\emph{Surface reactions.}
Reactions on the pore walls scaled by $\eps$, which give rise to volume
reaction terms in the limit, were first homogenized by Hornung and J\"ager
\cite{HornungJager1991} for adsorption and reaction of chemicals in porous
media; nonlinear monotone surface reactions were treated by Conca, D\'iaz and
Timofte \cite{ConcaDiazTimofte2003}. Piatnitski and Rybalko
\cite{PiatnitskiRybalko2011} considered rapidly oscillating Fourier-type
conditions of order one under a centering condition on the surface
nonlinearity, whereas the reaction here is of order $\eps$, need not be
centered, and yields a nonzero nonlinear volume reaction. Alternating
Robin--Neumann and dynamical boundary conditions were studied in
\cite{AmarAndreucciBellaveglia2017,Anguiano2020}, nonlinear boundary
reactions in evolving microdomains in \cite{GahnNeussRaduPop2021}, irregularly
perforated random domains in \cite{HeidaJahnelVu2025}, and nonlinear Robin
conditions on critically small holes, where a strange term appears, in
\cite{DiazPodolskiyShaposhnikova2025}.

\emph{Degenerate storage.}
Elliptic--parabolic equations arise naturally when the storage vanishes on
part of the domain \cite{AltLuckhaus1983}. Two-scale limits of models for
flow in partially fissured media were obtained by Clark and Showalter
\cite{ClarkShowalter1999} in the framework of implicit degenerate evolution
equations \cite{DiBenedettoShowalter1981}. In perforated
domains, a nonlinear storage $b(x/\eps,u)$, nondecreasing in $u$ with
$b(y,0)=0$, was considered in \cite{NandakumaranRajesh2002}, and a vanishing
heat capacity of order $\eps^q$ in \cite{DanielssonJohnsen2019}; here the
storage is linear and may vanish on a fixed part of each cell.

\emph{Degenerate diffusion weights.}
For degenerate weighted equations the notion of solution itself requires care,
since smooth functions need not be dense in the weighted Sobolev space
\cite{ChiadoPiatSerraCassano1994b,Zhikov1998}; relaxation of degenerate
variational integrals was studied in \cite{ChiadoPiatSerraCassano1994a}, and
periodic homogenization in
\cite{CardonePastukhovaPerugia2013,ZhikovPastukhova2008}, and random degenerate
monotone operators in \cite{EngstromPerssonPiatnitskiWall2006}. Integrability
conditions on the inverse of the coefficients, of the type used in $(H_2)$,
also appear in the stochastic theory
\cite{BellaFehrmanOtto2018,NeukammSchaffnerSchlomerkemper2017}.

\emph{Combined degeneracies.}
Both degeneracies occur in the work of Paronetto. In \cite{Paronetto2002}, the
same Muckenhoupt weight multiplies the time derivative and controls the
diffusion. In \cite{Paronetto2004}, periodic operators
$\mu(x)\partial_t-\Div(a(x,t)D)$ with $\mu\geq0$, possibly $\mu\equiv0$, are
homogenized in a domain without holes; the quadratic form of $a$ is comparable
to an $A_2$ weight, and a Muckenhoupt-type condition couples $\mu$ to this
weight. In a perforated domain, both degeneracies were treated in
\cite{DouanlaKenneBogningTetsadjio2026} for a quasilinear problem with a
Signorini condition on the exterior boundary.

The present setting differs from these works in three respects. The weights
satisfy only the directional inverse-integrability condition $(H_2)$, without
a Muckenhoupt condition, and the storage is subject only to
$\overline b>0$, without any compatibility condition with the weights. The
holes carry a monotone nonlinear Robin reaction of order $\eps$, whose limit is
a nonzero macroscopic volume reaction identified through strong convergence of
the scaled traces. Finally, we prove strong convergence of the weighted
gradients and uniform-in-time convergence of the energies and of the states in
the storage-weighted norm.

\subsection{Standing hypotheses}\label{sec:hypotheses}

All periodic coefficients below are identified with their periodic
extensions, and the subscript $\#$ denotes $Y$-periodicity: $C^\infty_\#(Y)$
is the space of smooth $Y$-periodic functions on $\R^N$; $L^p_\#(Y)$ and
$W^{1,p}_\#(Y)$ consist of the $Y$-periodic functions in
$L^p_{\rm loc}(\R^N)$ and $W^{1,p}_{\rm loc}(\R^N)$, respectively; and
$\mathscr D'_\#(Y)$ is the dual of $C^\infty_\#(Y)$.
We now state the assumptions on the diffusion, storage, reaction and data.

\noindent\textbf{(H$_2$).}\enspace
For some finite $s\geq N$, the diffusion weights satisfy
$w_i\in L^\infty_\#(Y)$ and
$w_i^{-s}\in L^1(Y)$ for every $i=1,\ldots,N$; in particular $w_i>0$
a.e. in $Y$. Write $\omega=(w_1,\ldots,w_N)$,
$A=\diag(w_1,\ldots,w_N)$ (so that
$(A^\eps)^{1/2}Dv=(\sqrt{w_i^\eps}D_iv)_{i=1}^N$), and put
\begin{equation}
 q=\frac{2s}{s+1}\in(1,2),\qquad
 q\geq\frac{2N}{N+1}>\frac{2N}{N+2}.
 \label{eq:q}
\end{equation}
The diffusion matrix is assumed diagonal.

\noindent\textbf{(H$_3$).}\enspace
The storage coefficient satisfies
\[
 b\in L^\infty_\#(Y),\qquad b\geq0\ \text{a.e. in }Y,\qquad
 \overline b:=\int_{\Ystar}b(y)\dd y>0.
\]
The values of $w_i$ and $b$ inside $\mathcal O$ are used on the unperforated
boundary strip, whereas their restrictions to $\Ystar$ determine the
cell averages. We set $\theta=|\Ystar|$. Since $|Y|=1$, neither
$\overline b$ nor the effective reaction is divided by $\theta$.

\noindent\textbf{(H$_4$).}\enspace
The surface reaction $h:S\times\R\to\R$ is measurable in $y$ and
extended periodically to
$\mathcal S=\bigcup_{k\in\mathbb Z^N}(k+S)$. There exists a constant
$L_h\geq0$, independent of $y$, such that, for a.e. $y\in S$ and
all $r,\widetilde r\in\R$,
\[
 h(y,0)=0,\qquad
 |h(y,r)-h(y,\widetilde r)|\leq L_h|r-\widetilde r|,\qquad
 [h(y,r)-h(y,\widetilde r)](r-\widetilde r)\geq0.
\]
Consequently, $|h(y,r)|\leq L_h|r|$ and $h(y,r)r\geq0$.

\noindent\textbf{(H$_5$).}\enspace
The data satisfy $f\in L^2(Q)$ and $u^{\rm in}\in L^2(\Om)$.

The assumptions play distinct roles: $(H_2)$ supplies weighted
gradient control, compact bulk embedding and scaled $L^2$ traces;
$\overline b>0$ in $(H_3)$ permits cell detection and produces the
positive effective storage; and $(H_4)$ provides microscopic
monotonicity and stability of the nonlinear surface reaction. Through the
monotonicity of $h^{\rm hom}$, it also yields uniqueness of the homogenized
problem, hence convergence of the whole family, and the sign condition
$h(y,r)r\geq0$ makes the cumulative surface dissipation nondecreasing in time,
which is used for the uniform-in-time convergence in
Section~\ref{sec:strong-corrector}.

\begin{remark}\label{rem:example}
For a spherical reference hole, let
\[
 w_i(y)=\min\{1,\dist(y,\mathcal S)^{\alpha_i}\},\qquad
 0<\alpha_i<1/N.
\]
Near $S$, the integral of $w_i^{-N}$ is controlled by a constant times
$\int_0^\delta r^{-N\alpha_i}\dd r$, which is finite. Thus $(H_2)$ holds
with $s=N$, and diffusion degenerates at the hole boundaries. If
$Y_b\subset\Ystar$ is measurable with $0<|Y_b|<|\Ystar|$, the periodic
extension of $b=\chi_{Y_b}$ satisfies $(H_3)$ and vanishes on a set of positive
measure in the material cell. The choice
$h(y,r)=c(y)\tanh r$, where $0\leq c\in L^\infty(S)$, satisfies $(H_4)$.
\end{remark}

\subsection{Main results and new ingredients}

\begin{maintheorem}\label{thm:A}
Assume $(H_1)$--$(H_5)$. For every $0<\eps<1$, problem \eqref{eq:micro}
has a unique solution $u_\eps$ in the sense of \eqref{eq:weak-micro}. As $\eps\to0$, the extended states
$P_\eps u_\eps$ converge strongly in $L^2(Q)$, for the whole family, to the
unique solution $u_0\in L^2(0,T;H^1_0(\Om))\cap H^1(0,T;H^{-1}(\Om))$ of
\[
 \overline b\,\partial_tu_0-\Div(A^{\rm hom}Du_0)+h^{\rm hom}(u_0)=\theta f
 \ \text{ in }Q,\qquad u_0=0\ \text{ on }\partial\Om\times(0,T),\qquad
 u_0(0)=u^{\rm in},
\]
where $\overline b=\int_{\Ystar}b$, $\theta=|\Ystar|$,
$h^{\rm hom}(r)=\int_Sh(y,r)\dd\sigma_y$, and the positive-definite matrix
$A^{\rm hom}$ is given by a weighted cell problem. Moreover, the dissipated
and stored energies converge uniformly on $[0,T]$, the weighted gradients
admit a strong corrector, and
\begin{equation*}
 \sup_{0\leq t\leq T}\int_{\Omeps}b^\eps(x)
       |U_\eps(x,t)-u_0(x,t)|^2\dd x\longrightarrow0
       \qquad(\eps\to0),
\end{equation*}
where $U_\eps(t)$ is the continuous representative of $u_\eps$ in the
storage-weighted space. These results are proved in
Theorems~\ref{thm:micro-solvability}, \ref{thm:hom-main}
and~\ref{thm:strong-corrector}.
\end{maintheorem}

The results identify the effective Robin reaction jointly with the
homogenized degenerate bulk dynamics, and establish uniform-in-time
energy and weighted-state convergence. The key ingredients are as follows.
\begin{enumerate}
\item A \emph{quantitative} cell-detection inequality
(Lemma~\ref{lem:cell-detection}), with rate $\eps^\beta$, transfers time
regularity of the mass to the state without a positive lower bound on $b$;
it sharpens, in the periodic setting, the qualitative detection inequality of
\cite{DouanlaKenneBogningTetsadjio2026}.
It yields time-translation estimates for the extended states
(Proposition~\ref{prop:time-translations}).
\item A scaled trace estimate (Lemma~\ref{lem:scaled-trace}) and strong
bulk compactness identify the nonlinear surface reaction directly
(Proposition~\ref{prop:surface-limit}). At the endpoint $s=N$, the cell
trace embedding into $L^2(S)$ is continuous but not generally compact;
strong scaled-trace convergence follows from the $\eps^2$-weighted
gradient remainder in the trace inequality.
\item The weighted two-scale gradient identification, based on an affine
closed-range argument
(Lemma~\ref{lem:twoscale-affine} and Theorem~\ref{thm:twoscale-weighted}),
recovers $H^1$ regularity of the limit and yields a positive-definite
effective diffusion matrix.
\item Energy comparison and terminal-mass duality
(Lemma~\ref{lem:terminal-storage}) give the strong weighted-gradient
corrector (Theorem~\ref{thm:strong-corrector}). Monotonicity of cumulative
dissipations, the energy identities and strong convergence of the mass
in $C([0,T];H^{-1}(\Om))$ yield uniform-in-time convergence of the
energies and of the states in the storage-weighted norm.
\end{enumerate}

\subsection{Organization of the paper}

Section~\ref{sec:setting} develops the weighted functional framework,
the scaled trace estimates and the nonlinear surface limit. Section~\ref{sec:micro}
treats the microscopic problem, and Section~\ref{sec:strong-compactness}
the detection and time-translation estimates and the compactness of the
states and of the weighted gradients.
Section~\ref{sec:homogenization} identifies the homogenized problem, and
Section~\ref{sec:strong-corrector} proves the uniform-in-time energy and
weighted-state convergence together with the strong gradient corrector.

\section{Functional framework and surface tools}\label{sec:setting}

The hypotheses of Section~\ref{sec:hypotheses} provide only weighted
energy bounds. We first develop the weighted energy spaces and extension
estimates, then establish the trace properties needed for the nonlinear Robin
reaction.

\subsection{Weighted functional setting}

For a fixed $\eps$, let $H^1_{\omega^\eps}(\Omeps)$ be the completion of
$C^\infty(\overline{\Omeps})$ for the norm
\begin{equation}
 \|v\|_{\omega^\eps}^2
 =\|v\|_{L^2(\Omeps)}^2+
 \sum_{i=1}^N\int_{\Omeps}w_i^\eps|D_i v|^2\dd x.
 \label{eq:energy-norm}
\end{equation}
Here smooth functions on a closure mean restrictions of functions smooth
in a neighbourhood of that closure; solutions are thus understood in the sense
of the completion, which avoids the Lavrentiev phenomenon that may occur for
degenerate weights
\cite{ChiadoPiatSerraCassano1994b,Zhikov1998,ZhikovPastukhova2008}. On the reference cell, denote by
$X_\omega(\Ystar)$ the analogous completion with norm
\[
 \|z\|_{\omega,\Ystar}^2
 =\|z\|_{L^2(\Ystar)}^2+
 \sum_{i=1}^N\int_{\Ystar}w_i|D_i z|^2\dd y.
\]

\begin{lemma}\label{lem:weighted-control}
Under $(H_1)$--$(H_2)$, the completions $H^1_{\omega^\eps}(\Omeps)$ and
$X_\omega(\Ystar)$ embed continuously into $W^{1,q}(\Omeps)$ and
$W^{1,q}(\Ystar)$, respectively. Moreover, for $i=1,\ldots,N$,
\begin{equation*}
 \|D_i v\|_{L^q(\Omeps)}
 \leq C\|\sqrt{w_i^\eps}D_i v\|_{L^2(\Omeps)}
 \qquad\forall\,v\in H^1_{\omega^\eps}(\Omeps),
\end{equation*}
with $C$ independent of $\eps$. On every fixed bounded Lipschitz domain $\mathcal G\subset\R^N$,
$W^{1,q}(\mathcal G)$ embeds compactly into $L^2(\mathcal G)$.
\end{lemma}

\begin{proof}
The identity $q/(2-q)=s$ and H\"older's inequality give
\[
 \|D_i v\|_{L^q(\Omeps)}
 \leq\left(\int_{\Omeps}(w_i^\eps)^{-s}\dd x\right)^{1/(2s)}
       \|\sqrt{w_i^\eps}D_i v\|_{L^2(\Omeps)}.
\]
Covering $\Om$ by at most $C_\Om\eps^{-N}$ periodic cells yields
$\int_{\Omeps}(w_i^\eps)^{-s}\leq C_\Om\int_Yw_i^{-s}$.
Since $q<2$, the $L^q(\Omeps)$ norm of the state is bounded by its
$L^2(\Omeps)$ norm. A Cauchy sequence for \eqref{eq:energy-norm} therefore
converges in $W^{1,q}(\Omeps)$. Its weighted derivatives converge in
$L^2(\Omeps)$, and extracting an a.e. convergent subsequence of the
ordinary derivatives shows that their limits are $\sqrt{w_i^\eps}D_iv$,
where $v$ is the $W^{1,q}(\Omeps)$ limit. Hence the natural map from the
completion into $W^{1,q}(\Omeps)$ is well defined, continuous and
injective: an element whose $W^{1,q}(\Omeps)$ limit vanishes has vanishing
weighted derivatives. We identify the completion with its image. The same
argument applies on $\Ystar$. Finally, \eqref{eq:q} gives the compact embedding
$W^{1,q}(\mathcal G)\Subset L^2(\mathcal G)$ for each fixed bounded Lipschitz domain $\mathcal G$
by the Sobolev--Rellich theorem; see \cite{AdamsFournier2003}.
\end{proof}

Define the Hilbert space
\begin{equation*}
 V_\eps=\{v\in H^1_{\omega^\eps}(\Omeps):
                 \gamma_0v=0\text{ on }\partial\Om\},
 \qquad \|v\|_{V_\eps}=\|v\|_{\omega^\eps},
\end{equation*}
where $\gamma_0:W^{1,q}(\Omeps)\to L^1(\partial\Om)$ is the trace on the
exterior boundary. Its closedness follows from
Lemma~\ref{lem:weighted-control} and continuity of the trace.

\begin{lemma}\label{lem:extension-poincare}
There are linear operators $P_\eps:V_\eps\to W^{1,q}_0(\Om)$ with
$P_\eps v=v$ on $\Omeps$ such that
\begin{equation}
 \|DP_\eps v\|_{L^q(\Om)}\leq C\|Dv\|_{L^q(\Omeps)},\qquad
 \|P_\eps v\|_{W^{1,q}(\Om)}\leq C\|v\|_{V_\eps}
 \qquad\forall\,v\in V_\eps.
 \label{eq:extension}
\end{equation}
In addition,
\begin{equation}
 \|v\|_{L^2(\Omeps)}^2
 \leq C_P\sum_{i=1}^N\int_{\Omeps}w_i^\eps|D_i v|^2\dd x
 \qquad\forall\,v\in V_\eps.
 \label{eq:poincare}
\end{equation}
All constants are independent of $\eps$.
\end{lemma}

\begin{proof}
Since $\overline{\mathcal O}\subset Y$, the set $\Ystar$ is a bounded,
connected open set with Lipschitz boundary $\partial Y\cup S$. By
\cite[Lemma~2.6]{AcerbiChiadoDalMasoPercivale1992}, applied with $A=\Ystar$,
$\omega=Y$ and $p=q$, there is a bounded linear extension operator
$E:W^{1,q}(\Ystar)\to W^{1,q}(Y)$ such that, for every $z\in W^{1,q}(\Ystar)$,
$Ez=z$ on $\Ystar$ (in particular near $\partial Y$) and
\[
 \|Ez\|_{L^q(Y)}\leq C\|z\|_{L^q(\Ystar)},\qquad
 \|D_yEz\|_{L^q(Y)}\leq C\|D_yz\|_{L^q(\Ystar)}.
\]
The gradient bound comes from applying a fixed extension to $z-z_{\Ystar}$,
where $z_{\Ystar}$ is the mean of $z$ over $\Ystar$, and using the
Poincar\'e--Wirtinger inequality on $\Ystar$; this is the construction of
\cite[Lemma~3]{CioranescuSaintJeanPaulin1979}, given there for $q=2$ and
smooth holes.

For $v\in V_\eps$, define $P_\eps v$ on each complete cell $\eps(k+Y)$,
$k\in I_\eps$, by $P_\eps v(x)=(Ez_k)(x/\eps-k)$ with $z_k(y)=v(\eps(k+y))$, and
set $P_\eps v=v$ on the unperforated strip $\Sigma_\eps$, as in the proof of
\cite[Theorem~4]{CioranescuSaintJeanPaulin1979}. Since $Ez_k=z_k$ near
$\partial Y$, the cellwise definitions agree with $v$ on all cell faces, so
$P_\eps v\in W^{1,q}(\Om)$ and $P_\eps v=v$ on $\Omeps$. Under $x=\eps(k+y)$
both gradient integrals acquire the factor $\eps^{N-q}$, so that
$\int_{\eps(k+Y)}|DP_\eps v|^q\dd x\leq C\int_{\eps(k+\Ystar)}|Dv|^q\dd x$.
Summing over $k\in I_\eps$ gives the first estimate in \eqref{eq:extension},
and the rescaled $L^q$ bound gives
$\|P_\eps v\|_{L^q(\Om)}\leq C\|v\|_{L^q(\Omeps)}$.
The construction affects only holes strictly inside $\Om$, hence the
exterior trace remains $\gamma_0v=0$. The weighted H\"older estimate
then proves the second bound in \eqref{eq:extension}.
Poincar\'e's inequality on $\Om$, followed by the Sobolev embedding,
now gives
\[
 \|v\|_{L^2(\Omeps)}\leq\|P_\eps v\|_{L^2(\Om)}
 \leq C\|DP_\eps v\|_{L^q(\Om)}
 \leq C\Bigl(\sum_{i=1}^N\int_{\Omeps}w_i^\eps|D_i v|^2\Bigr)^{1/2}.
\]
This proves \eqref{eq:poincare}.
\end{proof}

\begin{lemma}\label{lem:periodic-averaging}
Under $(H_1)$, the strip $\Sigma_\eps$ is covered, up to a null set, by at
most $C_\Om\eps^{1-N}$ cubes $\eps(k+Y)$, $k\in\mathbb Z^N\setminus I_\eps$;
in particular $|\Sigma_\eps|\leq C_\Om\eps$. For every $1\leq p<\infty$ and
$a\in L^p_\#(Y)$,
\begin{equation}
 \lim_{\eps\to0}\int_\Om|a(x/\eps)|^p\dd x
     =|\Om|\int_Y|a(y)|^p\dd y,\qquad
 \int_{\Sigma_\eps}|a(x/\eps)|^p\dd x
     \leq C\eps\|a\|_{L^p(Y)}^p,
 \label{eq:periodic-cell-average}
\end{equation}
and $\int_\Om|a(x/\eps)|^p\dd x\leq C\|a\|_{L^p(Y)}^p$ for all
$0<\eps<1$. If $a\in L^1_\#(Y)$, then
\begin{equation}
 \int_Q a(x/\eps)\phi(x,t)\dd x\dd t\to\int_Ya(y)\dd y\int_Q\phi\dd x\dd t
 \qquad\forall\,\phi\in C(\overline Q).
 \label{eq:periodic-L1-average}
\end{equation}
If $a\in L^\infty_\#(Y)$, then
\begin{equation}
 a(x/\eps)\weakto^*\int_Ya(y)\dd y,\qquad
 \chi_{\Omeps}a(x/\eps)\weakto^*\int_{\Ystar}a(y)\dd y
 \quad\text{in }L^\infty(\Om).
 \label{eq:periodic-material-average}
\end{equation}
Finally, let $1<p<\infty$, $p'=p/(p-1)$ and $a\in L^{p'}_\#(Y)$. If $z_\eps$
is bounded in $L^p(Q)$ and converges weakly two-scale to $z$ in $L^p(Q)$,
then
\begin{equation}
 \int_Q z_\eps(x,t)a(x/\eps)\Theta(x,t,x/\eps)\dd x\dd t
 \to\int_Q\int_Yz(x,t,y)a(y)\Theta(x,t,y)\dd y\dd x\dd t
 \qquad\forall\,\Theta\in C_c^\infty\bigl(Q;C^\infty_\#(Y)\bigr).
 \label{eq:periodic-multiplier}
\end{equation}
Consequently, if $z_\eps a(x/\eps)$ is bounded in $L^r(Q)$ for some
$1<r<\infty$, then $z_\eps a(x/\eps)$ converges weakly two-scale to $za$
in $L^r(Q)$. Moreover, if $(\zeta_\eps)$ is bounded in $L^p(Q)$, then
\begin{equation}
 \|\mathbf 1_{\Sigma_\eps}\zeta_\eps\|_{L^1(Q)}
 \leq(T|\Sigma_\eps|)^{1/p'}\|\zeta_\eps\|_{L^p(Q)}\leq C\eps^{1/p'},
 \label{eq:strip-error}
\end{equation}
so that $\mathbf 1_{\Sigma_\eps}\zeta_\eps$ has zero weak two-scale limit in
$L^p(Q)$. The assertions hold componentwise for vector- and matrix-valued
periodic functions.
\end{lemma}

\begin{proof}
Outside cell interfaces, each $x\in\Sigma_\eps$ lies in a cube
$\eps(k+Y)$ with $k\notin I_\eps$ that meets $\Om$. Such a cube is not
contained in $\Om$, hence meets $\partial\Om$, and lies in the
$2\sqrt N\,\eps$-neighbourhood of $\partial\Om$. The bounded Lipschitz
boundary is covered by finitely many graph charts; in each chart its
two-sided $\eta$-neighbourhood lies in a slab of thickness at most $C\eta$
about the graph, whence
$|\{x\in\R^N:\dist(x,\partial\Om)<\eta\}|\leq C_\Om\eta$ for $0<\eta<\eta_0$. The cubes
have disjoint interiors and volume $\eps^N$, so their number $n_\eps$
satisfies $n_\eps\eps^N\leq C_\Om\eps$, that is, $n_\eps\leq C_\Om\eps^{1-N}$
(with an adjusted constant for $\eps\geq\eta_0/(2\sqrt N)$). This proves the
first assertion. Periodicity and a change of variables on those cubes give
\[
 \int_{\Sigma_\eps}|a(x/\eps)|^p\dd x
 \leq C\eps^{1-N}\eps^N\int_Y|a(y)|^p\dd y,
\]
which is the second assertion in \eqref{eq:periodic-cell-average}.
The complete cells contribute $\eps^N\#I_\eps\int_Y|a|^p$, and
$\eps^N\#I_\eps\to|\Om|$ because $|\Sigma_\eps|=O(\eps)$; this gives the
limit and the uniform bound.

For a smooth periodic $a$, decomposition into complete cells and the
uniform continuity of $\phi\in C(\overline Q)$ give \eqref{eq:periodic-L1-average};
for $a\in L^1_\#(Y)$, approximation in $L^1(Y)$ by smooth periodic functions
and the uniform bound with $p=1$ remove the error, uniformly in $\eps$.
Replacing $a$ by $\chi_{\Ystar}a$ gives the material-cell counterpart. For
$a\in L^\infty_\#(Y)$, density of continuous functions in $L^1(\Om)$ then
yields \eqref{eq:periodic-material-average}. On the incomplete cells
$\chi_{\Omeps}$ equals one rather than $\chi_{\Ystar}(x/\eps)$; the
difference is bounded and supported in $\Sigma_\eps$, so it disappears in
$L^1(\Om)$.

Finally choose $a_m\in C^\infty_\#(Y)$ with $a_m\to a$ in $L^{p'}(Y)$. By the
uniform bound, for $0<\eps<1$,
\[
 \|(a-a_m)(x/\eps)\Theta(x,t,x/\eps)\|_{L^{p'}(Q)}
 \leq C\|\Theta\|_{L^\infty(Q\times Y)}\|a-a_m\|_{L^{p'}(Y)}.
\]
For fixed $m$, $a_m\Theta$ is an admissible test, and the limit follows
from the definition of weak two-scale convergence. The corresponding
limit-side difference is bounded by
$C\|z\|_{L^p(Q\times Y)}\|\Theta\|_{L^\infty(Q\times Y)}\|a-a_m\|_{L^{p'}(Y)}$.
Passing first $\eps\to0$ and then $m\to\infty$ proves
\eqref{eq:periodic-multiplier}. If $z_\eps a(x/\eps)$ is bounded in
$L^r(Q)$, each subsequence has a weak two-scale limit in $L^r(Q\times Y)$,
which \eqref{eq:periodic-multiplier} identifies as $za$. Finally,
\eqref{eq:strip-error} is H\"older's inequality together with
$|\Sigma_\eps|\leq C_\Om\eps$; the $L^1(Q)$ convergence identifies the
two-scale limit of the bounded sequence $\mathbf 1_{\Sigma_\eps}\zeta_\eps$
against smooth tests, and density gives the general case.
\end{proof}

\subsection{Scaled traces and the boundary operator}

Since the total area of the holes is of order $\eps^{-1}$, traces on
$S_\eps$ are measured with the weight $\eps$. The following estimate is
uniform in $\eps$ and involves only weighted gradients.

\begin{lemma}\label{lem:scaled-trace}
Under $(H_1)$--$(H_2)$, the trace $\gamma_S:X_\omega(\Ystar)\to L^2(S)$
is continuous. For every $v\in H^1_{\omega^\eps}(\Omeps)$,
\begin{equation}
 \eps\int_{S_\eps}|\gamma_\eps v|^2\dd\sigma_x
 \leq C\left(\int_{\Omeps}|v|^2\dd x
       +\eps^2\sum_{i=1}^N\int_{\Omeps}w_i^\eps|D_i v|^2\dd x\right),
 \label{eq:scaled-trace}
\end{equation}
where $\gamma_\eps$ is the trace on the holes and $C$ is independent of
$\eps$.
\end{lemma}

\begin{proof}
Since $1<q<N$, the continuous Sobolev trace embedding has exponent
$q_\partial=q(N-1)/(N-q)$. The condition $q\geq2N/(N+1)$ gives
$q_\partial\geq2$. The trace theorem and the cell version of
Lemma~\ref{lem:weighted-control} yield, for every $z\in W^{1,q}(\Ystar)$
with finite weighted norm,
\begin{equation}
 \|\gamma_S z\|_{L^2(S)}^2
 \leq C\|z\|_{W^{1,q}(\Ystar)}^2
 \leq C\left(\|z\|_{L^2(\Ystar)}^2+
                  \sum_{i=1}^N\int_{\Ystar}w_i|D_i z|^2\dd y\right).
 \label{eq:cell-trace}
\end{equation}
The trace extends to the weighted completion by
\eqref{eq:cell-trace}.

For $k\in I_\eps$, apply \eqref{eq:cell-trace} to
$z_k(y)=v(\eps(k+y))$. Changing variables gives
\[
 \eps\|\gamma_\eps v\|_{L^2(\eps(k+S))}^2
 \leq C\bigl(\|v\|_{L^2(\eps(k+\Ystar))}^2
    +\eps^2\|(A^\eps)^{1/2}Dv\|_{L^2(\eps(k+\Ystar))^N}^2\bigr).
\]
Summing over the disjoint complete cells proves \eqref{eq:scaled-trace}
for smooth functions. Approximation in the weighted norm proves it for
the completion. The unperforated strip contributes no surface term.
\end{proof}

The weak form of \eqref{eq:micro} involves a mass operator and an energy
operator: for $v,z\in V_\eps$, we define
\begin{align}
 \langle \mathcal B_\eps v,z\rangle
   &=\int_{\Omeps}b^\eps vz\dd x,\nonumber\\
 \langle \mathcal A_\eps v,z\rangle
   &=\sum_{i=1}^N\int_{\Omeps}w_i^\eps D_i vD_i z\dd x
     +\eps\int_{S_\eps}
          h(x/\eps,\gamma_\eps v)\gamma_\eps z\dd\sigma_x.
 \label{eq:operators}
\end{align}
The coefficient $b$ may vanish, so $\mathcal B_\eps$ need not be injective.

\begin{proposition}\label{prop:operator}
Under $(H_1)$, $(H_2)$ and $(H_4)$, $\mathcal A_\eps:V_\eps\to V_\eps'$ is globally
Lipschitz and strongly monotone, with constants independent of $\eps$:
\begin{equation}
 \|\mathcal A_\eps v-\mathcal A_\eps z\|_{V_\eps'}\leq C\|v-z\|_{V_\eps},\qquad
 \langle \mathcal A_\eps v-\mathcal A_\eps z,v-z\rangle
       \geq \frac{\|v-z\|_{V_\eps}^2}{1+C_P}
 \qquad\forall\,v,z\in V_\eps.
 \label{eq:operator-monotone}
\end{equation}
Moreover $\mathcal A_\eps0=0$.
\end{proposition}

\begin{proof}
For $\zeta\in V_\eps$, the difference of the surface terms is bounded by
\[
 L_h\left(\eps\int_{S_\eps}|\gamma_\eps(v-z)|^2\right)^{1/2}
       \left(\eps\int_{S_\eps}|\gamma_\eps\zeta|^2\right)^{1/2}
 \leq C\|v-z\|_{V_\eps}\|\zeta\|_{V_\eps},
\]
by \eqref{eq:scaled-trace} and $\eps<1$. Cauchy--Schwarz gives the same
bound for the diffusion term. Monotonicity of $h$ gives
\[
 \langle \mathcal A_\eps v-\mathcal A_\eps z,v-z\rangle
 \geq\sum_{i=1}^N\int_{\Omeps}w_i^\eps|D_i(v-z)|^2,
\]
and \eqref{eq:poincare} proves \eqref{eq:operator-monotone}.
The normalization $h(y,0)=0$ gives $\mathcal A_\eps0=0$.
\end{proof}

\subsection{Strong traces and the nonlinear surface limit}

Lemma~\ref{lem:scaled-trace} bounds the traces. To pass to the limit in
the reaction, the next result derives their strong convergence from bulk
convergence, using Riemann sums instead of two-scale convergence on surfaces
\cite{AllaireDamlamianHornung1996,NeussRadu1996}. We suppress the trace
symbols inside surface integrals.

\begin{proposition}\label{prop:surface-limit}
Assume $(H_1)$--$(H_2)$ and $(H_4)$. Let
$u_\eps\in L^2(0,T;H^1_{\omega^\eps}(\Omeps))$ and
$u\in L^2(0,T;H^1(\Om))$ satisfy
\begin{equation}
 \int_{\Qeps}|u_\eps-u|^2\dd x\dd t\to0,
 \qquad
 \sup_{0<\eps<1}\sum_{i=1}^N\int_{\Qeps}w_i^\eps|D_i u_\eps|^2\dd x\dd t<\infty.
 \label{eq:surface-premises}
\end{equation}
Then
\begin{align}
 \eps\int_0^T\int_{S_\eps}|u_\eps-u|^2\dd\sigma_x\dd t&\to0,
 \label{eq:strong-surface}\\
 \eps\int_0^T\int_{S_\eps}
 |h(x/\eps,u_\eps)-h(x/\eps,u)|^2\dd\sigma_x\dd t&\to0.
 \label{eq:strong-reaction}
\end{align}
Define
\begin{equation}
 h^{\rm hom}(r)=\int_S h(y,r)\dd\sigma_y\qquad\forall\,r\in\R.
 \label{eq:effective-reaction}
\end{equation}
Then
\begin{equation}
 \eps\int_0^T\int_{S_\eps}
       h(x/\eps,u_\eps)\varphi\dd\sigma_x\dd t
 \to\int_Q h^{\rm hom}(u)\varphi\dd x\dd t
 \qquad\forall\,\varphi\in L^2(0,T;H^1(\Om)).
 \label{eq:reaction-limit}
\end{equation}
The surface dissipation also converges:
\begin{equation}
 \eps\int_0^T\int_{S_\eps}
       h(x/\eps,u_\eps)u_\eps\dd\sigma_x\dd t
 \to\int_Q h^{\rm hom}(u)u\dd x\dd t.
 \label{eq:dissipation-limit}
\end{equation}
\end{proposition}

\begin{proof}
Since the weights are bounded, the restriction of $u$ to $\Qeps$ belongs
to $L^2(0,T;H^1_{\omega^\eps}(\Omeps))$. Apply \eqref{eq:scaled-trace} to $u_\eps-u$ and integrate
in time:
\[
 \eps\|u_\eps-u\|_{L^2(S_\eps\times(0,T))}^2
 \leq C\bigl(\|u_\eps-u\|_{L^2(\Qeps)}^2
 +\eps^2\|(A^\eps)^{1/2}Du_\eps\|_{L^2(\Qeps)^N}^2
 +\eps^2\|Du\|_{L^2(Q)^N}^2\bigr).
\]
Every term tends to zero, which proves \eqref{eq:strong-surface}.
The Lipschitz bound on $h$ immediately gives
\eqref{eq:strong-reaction}.

Taking all weights equal to one in Lemma~\ref{lem:scaled-trace}, for
$z\in L^2(0,T;H^1(\Om))$ we obtain
\begin{equation}
 \eps\int_0^T\int_{S_\eps}|z|^2\dd\sigma_x\dd t
 \leq C\bigl(\|z\|_{L^2(Q)}^2+
                   \eps^2\|Dz\|_{L^2(Q)}^2\bigr).
 \label{eq:macro-trace}
\end{equation}
By \eqref{eq:strong-reaction}, \eqref{eq:macro-trace} and
Cauchy--Schwarz, replacing $u_\eps$ by $u$ in the left-hand side of
\eqref{eq:reaction-limit} changes it by a quantity tending to zero.

Suppose first that $u$ and $\varphi$ are smooth on a neighbourhood of
$\overline\Om\times[0,T]$. Changing variables on each hole gives
\[
 \eps^N\sum_{k\in I_\eps}\int_0^T\int_S
 h(y,u(\eps k+\eps y,t))\varphi(\eps k+\eps y,t)
                       \dd\sigma_y\dd t.
\]
Replace the two macroscopic arguments by their values at $(\eps k,t)$.
Uniform continuity, the Lipschitz bound on $h$, and
$\eps^N\#I_\eps\leq|\Om|$ show that the error tends to zero.
The resulting Riemann sum converges to
$\int_Q h^{\rm hom}(u)\varphi$; the omitted strip has measure $O(\eps)$.

Set $\mathcal X=L^2(0,T;H^1(\Om))$. For $u,\varphi\in \mathcal X$, choose smooth
$u_\ell,\varphi_\ell$ converging to them in $\mathcal X$. By
\eqref{eq:macro-trace} and the Lipschitz bound, uniformly in
$0<\eps<1$,
\[
 \left|\eps\int_0^T\int_{S_\eps}
 [h(x/\eps,u)\varphi-h(x/\eps,u_\ell)\varphi_\ell]\right|
 \leq C\bigl(\|u-u_\ell\|_{\mathcal X}\|\varphi\|_{\mathcal X}
 +\|u_\ell\|_{\mathcal X}\|\varphi-\varphi_\ell\|_{\mathcal X}\bigr).
\]
The same approximation works in the bulk, since $h^{\rm hom}$ is
Lipschitz with constant $L_h|S|$ and vanishes at zero. Passing first
$\eps\to0$ and then $\ell\to\infty$ proves \eqref{eq:reaction-limit}.

Finally, \eqref{eq:scaled-trace}, \eqref{eq:surface-premises} and
$|h(y,r)|\leq L_h|r|$ show that
$\eps\int_0^T\int_{S_\eps}|h(x/\eps,u_\eps)|^2$ is bounded.
Together with \eqref{eq:strong-surface}, this implies
\[
 \eps\int_0^T\int_{S_\eps}
 h(x/\eps,u_\eps)(u_\eps-u)\dd\sigma_x\dd t\to0.
\]
Apply \eqref{eq:reaction-limit} with $\varphi=u$ to conclude
\eqref{eq:dissipation-limit}.
\end{proof}

\section{The microscopic problem and solvability}\label{sec:micro}

We establish solvability by regularizing the mass operator, within the
framework of implicit degenerate evolution equations
\cite{DiBenedettoShowalter1981}. This section requires only $b\in L^\infty_\#(Y)$ with $b\geq0$; the condition
$\overline b>0$ is reserved for the homogenization problem.

The operators in \eqref{eq:operators} specify the weak formulation of \eqref{eq:micro}: seek
\begin{equation}
 u_\eps\in L^2(0,T;V_\eps),\qquad
 \mathcal B_\eps u_\eps\in H^1(0,T;V_\eps')
 \label{eq:solution-space}
\end{equation}
such that
\begin{equation}
 (\mathcal B_\eps u_\eps)' + \mathcal A_\eps u_\eps=F_\eps
 \quad\text{in }L^2(0,T;V_\eps'),\qquad
 (\mathcal B_\eps u_\eps)(0)=g_\eps^{\rm in},
 \label{eq:weak-micro}
\end{equation}
where $\langle F_\eps(t),z\rangle=\int_{\Omeps}f(x,t)z(x)\dd x$ and
$\langle g_\eps^{\rm in},z\rangle=\int_{\Omeps}b^\eps u^{\rm in}z\dd x$; that is,
for a.e. $t\in(0,T)$,
\begin{equation*}
 \langle(\mathcal B_\eps u_\eps)'(t),z\rangle
 +\sum_{i=1}^N\int_{\Omeps}w_i^\eps D_iu_\eps D_iz\dd x
 +\eps\int_{S_\eps}h(x/\eps,u_\eps)z\dd\sigma_x
 =\int_{\Omeps}fz\dd x
 \qquad\forall\,z\in V_\eps.
\end{equation*}
The initial mass is defined by the continuous $V_\eps'$ representative in
\eqref{eq:solution-space}.

\subsection{The mass space and the energy identity}

For a fixed $\eps$, let
\begin{equation*}
 \mathcal K_\eps=\ker \mathcal B_\eps
 =\left\{v\in V_\eps:\int_{\Omeps}b^\eps|v|^2\dd x=0\right\},
 \qquad
 H_{b,\eps}=\overline{V_\eps/\mathcal K_\eps}^{\,\|\cdot\|_{b,\eps}},
\end{equation*}
where $\|[v]\|_{b,\eps}^2=\int_{\Omeps}b^\eps|v|^2\dd x$.
The natural map $J_\eps:V_\eps\to H_{b,\eps}$ is continuous and has dense
range. Its adjoint is understood as a map into $V_\eps'$:
\[
 \langle J_\eps^*\mu,v\rangle=(\mu,J_\eps v)_{H_{b,\eps}}.
\]
Thus $\mathcal B_\eps=J_\eps^*J_\eps$, and $J_\eps^*$ is injective. The kernel
of $J_\eps$ is exactly $\mathcal K_\eps$. Since $C_c^\infty(\Omeps)\subset V_\eps$
is dense in $L^2(\Omeps,b^\eps\dd x)$, the space $H_{b,\eps}$ is isometric to
$L^2(\Omeps,b^\eps\dd x)$, where functions are identified
$b^\eps\dd x$-almost everywhere; we use this identification throughout. In
particular, $u^{\rm in}|_{\Omeps}\in H_{b,\eps}$, and the initial mass in
\eqref{eq:weak-micro} is $g_\eps^{\rm in}=J_\eps^*u^{\rm in}$.

\begin{lemma}\label{lem:mass-energy}
Let $v\in L^2(0,T;V_\eps)$ and
$(\mathcal B_\eps v)'\in L^2(0,T;V_\eps')$. Then $J_\eps v$ has a representative
in $C([0,T];H_{b,\eps})$, and
\begin{equation}
 \frac12\|J_\eps v(t_1)\|_{H_{b,\eps}}^2
 -\frac12\|J_\eps v(t_0)\|_{H_{b,\eps}}^2
 =\int_{t_0}^{t_1}\langle(\mathcal B_\eps v)'(t),v(t)\rangle\dd t
 \qquad\forall\,0\leq t_0\leq t_1\leq T.
 \label{eq:mass-energy}
\end{equation}
Moreover,
\begin{equation}
 (\mathcal B_\eps v)(t)=J_\eps^*(J_\eps v)(t)
 \quad\text{for every }t\in[0,T],\qquad
 (\mathcal B_\eps v)(0)=J_\eps^*\mu\ \Longleftrightarrow\ (J_\eps v)(0)=\mu
 \quad\forall\,\mu\in H_{b,\eps}.
 \label{eq:mass-initial-equivalence}
\end{equation}
Here $J_\eps v(t)$ denotes the continuous mass-space representative.
\end{lemma}

\begin{proof}
We follow \cite[Proposition~III.1.2 and Section~III.6]{Showalter1997}.
The kernel $\mathcal K_\eps$ is closed in $V_\eps$. Let
$\widetilde V_\eps=\mathcal K_\eps^\perp$, with orthogonality taken in $V_\eps$, and let
$\Pi_\eps$ be the orthogonal projection onto $\widetilde V_\eps$.
The restriction of $J_\eps$ to $\widetilde V_\eps$ is injective, continuous and has
dense range in $H_{b,\eps}$, because $J_\eps\Pi_\eps=J_\eps$.
It therefore defines the Hilbert triple
$\widetilde V_\eps\hookrightarrow H_{b,\eps}\hookrightarrow \widetilde V_\eps'$.

Put $\varpi=(\mathcal B_\eps v)'$. Since $\mathcal B_\eps v$ annihilates $\mathcal K_\eps$, so
does $\varpi$ for almost every time. The Gelfand embedding
$H_{b,\eps}\hookrightarrow \widetilde V_\eps'$ sends $J_\eps x$ to the
functional $y\mapsto(J_\eps x,J_\eps y)_{H_{b,\eps}}$, namely
$(\mathcal B_\eps x)|_{\widetilde V_\eps}$. Thus
\[
 -\int_0^T (J_\eps v,J_\eps y)_{H_{b,\eps}}\phi'(t)\dd t
 =\int_0^T\langle \varpi(t),y\rangle\phi(t)\dd t
 \qquad\forall\,y\in \widetilde V_\eps,\ \forall\,\phi\in C_c^\infty(0,T).
\]
Consequently the time derivative of the canonical image of
$\Pi_\eps v\in L^2(0,T;\widetilde V_\eps)$ is exactly
$\varpi|_{\widetilde V_\eps}\in L^2(0,T;\widetilde V_\eps')$.
The Hilbert-triple energy identity
\cite[Chapter~1, Theorem~3.1]{LionsMagenes1972},
\cite[Proposition~III.1.2]{Showalter1997} therefore gives continuity of
$J_\eps\Pi_\eps v=J_\eps v$ in $H_{b,\eps}$ and
\[
 \frac12\frac{\dd}{\dd t}\|J_\eps v\|_{H_{b,\eps}}^2
 =\langle \varpi,\Pi_\eps v\rangle=\langle \varpi,v\rangle
 \quad\text{a.e. in }(0,T).
\]
This proves \eqref{eq:mass-energy}. The function $\mathcal B_\eps v$ belongs to
$H^1(0,T;V_\eps')$. Its continuous representative agrees a.e., and hence
everywhere, with $J_\eps^*(J_\eps v)$. Injectivity of $J_\eps^*$ proves
the equivalence in \eqref{eq:mass-initial-equivalence}.
\end{proof}

\subsection{Existence, uniqueness and uniform estimates}

Lemma~\ref{lem:mass-energy} provides the energy identity needed both for
existence and for estimates uniform in $\eps$. These estimates are stated for
the following quantities, defined on $\Om$:
\[
 \bar u_\eps=P_\eps u_\eps,\qquad
 \rho_\eps=\chi_{\Omeps}b^\eps,\qquad
 m_\eps=\rho_\eps\bar u_\eps=\chi_{\Omeps}b^\eps u_\eps;
\]
the extension by zero is applied to the mass, not to the state.

\begin{theorem}\label{thm:micro-solvability}
Assume $(H_1)$, $(H_2)$ and $(H_4)$, let
$b\in L^\infty_\#(Y)$ be nonnegative, and let $f\in L^2(Q)$ and
$u^{\rm in}\in L^2(\Om)$. For every $0<\eps<1$, problem
\eqref{eq:weak-micro} has exactly one solution in the class
\eqref{eq:solution-space}. It satisfies $(J_\eps u_\eps)(0)=u^{\rm in}$ and,
for every $t\in[0,T]$,
\begin{align}
 \frac12\|J_\eps u_\eps(t)\|_{H_{b,\eps}}^2
 &+\int_0^t\sum_{i=1}^N\int_{\Omeps}w_i^\eps|D_i u_\eps|^2\dd x\dd r
 \nonumber\\
 &+\eps\int_0^t\int_{S_\eps}
       h(x/\eps,u_\eps)u_\eps\dd\sigma_x\dd r
 =\frac12\int_{\Omeps}b^\eps|u^{\rm in}|^2\dd x
 +\int_0^t\int_{\Omeps}fu_\eps\dd x\dd r.
 \label{eq:micro-energy}
\end{align}
There is a constant $C$, independent of $\eps$, $f$ and $u^{\rm in}$, such
that, with $M_0^2=\|f\|_{L^2(Q)}^2+\|u^{\rm in}\|_{L^2(\Om)}^2$,
\begin{equation}
 \|u_\eps\|_{L^2(0,T;V_\eps)}^2
 +\sup_{0\leq t\leq T}\|J_\eps u_\eps(t)\|_{H_{b,\eps}}^2
 +\|(\mathcal B_\eps u_\eps)'\|_{L^2(0,T;V_\eps')}^2
 \leq CM_0^2.
 \label{eq:micro-bound}
\end{equation}
Moreover,
\begin{align}
 &\|\bar u_\eps\|_{L^2(0,T;W^{1,q}_0(\Om))}^2
   +\sum_{i=1}^N\|\chi_{\Omeps}w_i^\eps D_i u_\eps\|_{L^2(Q)}^2
   \leq CM_0^2,\label{eq:uniform-bulk}\\
 &\eps\int_0^T\int_{S_\eps}
        \bigl(|u_\eps|^2+|h(x/\eps,u_\eps)|^2\bigr)
           \dd\sigma_x\dd t
  +\|m_\eps\|_{H^1(0,T;H^{-1}(\Om))}^2\leq CM_0^2.
                \label{eq:uniform-mass}
\end{align}
\end{theorem}

\begin{proof}
Fix $\eps$ and abbreviate $V=V_\eps$, $H=H_{b,\eps}$,
$\mathcal B=\mathcal B_\eps$, $J=J_\eps$, $\mathcal A=\mathcal A_\eps$ and $F=F_\eps$; note that
$\|F\|_{L^2(0,T;V')}\leq\|f\|_{L^2(Q)}$. By Proposition~\ref{prop:operator},
$\mathcal A$ is Lipschitz with constant $C_{\mathcal A}$ and strongly monotone with constant
$c_*=(1+C_P)^{-1}$, both independent of $\eps$, and $\mathcal A0=0$. For existence
of solutions of $(\mathcal B u)'+\mathcal A u=F$ with $\mathcal B$ symmetric
and nonnegative we follow \cite{DiBenedettoShowalter1981} and
\cite[Section~III.6]{Showalter1997}; the argument also yields the energy
identity. Let
$\mathcal I:V\to V'$ be the Riesz isomorphism. For $0<\delta\leq1$, the
operator $\mathcal B_\delta=\mathcal B+\delta\mathcal I$ is an isomorphism, and Picard
iteration gives a unique $u_\delta\in H^1(0,T;V)$ with
$\mathcal B_\delta u_\delta'+\mathcal Au_\delta=F$ and
$u_\delta(0)=\mathcal B_\delta^{-1}g_\eps^{\rm in}$. Write
$e_\delta(t)=\|Ju_\delta(t)\|_H^2+\delta\|u_\delta(t)\|_V^2$. Pairing the
initial condition with $u_\delta(0)$ gives
\[
 e_\delta(0)=\langle g_\eps^{\rm in},u_\delta(0)\rangle
 =(u^{\rm in},Ju_\delta(0))_H\leq\|u^{\rm in}\|_H\,e_\delta(0)^{1/2},
\]
so that $e_\delta(0)\leq\|u^{\rm in}\|_H^2$. Pairing the equation with
$u_\delta$ and Young's inequality then give
\[
 \sup_{0\leq t\leq T}e_\delta(t)
 +\|u_\delta\|_{L^2(0,T;V)}^2\leq C\bigl(\|F\|_{L^2(0,T;V')}^2+\|u^{\rm in}\|_H^2\bigr),
\]
so that $m_\delta=\mathcal B_\delta u_\delta$ is bounded in $H^1(0,T;V')$ and
$\delta\mathcal Iu_\delta\to0$ in $C([0,T];V')$. After extraction,
$u_\delta\weakto u$ in $L^2(0,T;V)$, $\mathcal Au_\delta\weakto\mathfrak a$ in
$L^2(0,T;V')$ and $m_\delta\weakto \mathcal Bu$ in $H^1(0,T;V')$, so that
$(\mathcal Bu)'=F-\mathfrak a$ and $(\mathcal Bu)(0)=g_\eps^{\rm in}$; by
\eqref{eq:mass-initial-equivalence} and injectivity of $J^*$,
$Ju(0)=u^{\rm in}$. In particular,
Lemma~\ref{lem:mass-energy} supplies a well-defined terminal value
$Ju(T)\in H$ and identifies $\mathcal Bu(T)=J^*Ju(T)$ in $V'$.
The regularized energy bound makes $Ju_\delta(T)$ bounded in $H$.
If $\lambda$ is any weak cluster point in $H$, continuity of the endpoint
trace from $H^1(0,T;V')$ to $V'$ and the estimate
$\|\delta\mathcal Iu_\delta(T)\|_{V'}\leq C\sqrt\delta$
give
\[
 J^*Ju_\delta(T)=m_\delta(T)-\delta\mathcal Iu_\delta(T)
 \rightharpoonup \mathcal Bu(T)=J^*Ju(T)
 \quad\text{in }V'.
\]
The injectivity of $J^*$ implies $\lambda=Ju(T)$; therefore
$Ju_\delta(T)\rightharpoonup Ju(T)$ in $H$. In the same way,
$J^*Ju_\delta(0)=g_\eps^{\rm in}-\delta\mathcal Iu_\delta(0)\to g_\eps^{\rm in}=J^*u^{\rm in}$
in $V'$, so that $Ju_\delta(0)\rightharpoonup u^{\rm in}$ in $H$ and
$e_\delta(0)=(u^{\rm in},Ju_\delta(0))_H\to\|u^{\rm in}\|_H^2$. Since
$\int_0^T\langle\mathcal Au_\delta,u_\delta\rangle
=\frac12e_\delta(0)-\frac12e_\delta(T)+\int_0^T\langle F,u_\delta\rangle$,
weak lower semicontinuity and the mass-energy identity for $u$ give
\[
 \limsup_{\delta\to0}\int_0^T\langle \mathcal Au_\delta,u_\delta\rangle
 \leq\frac12\|u^{\rm in}\|_H^2-\frac12\|Ju(T)\|_H^2+\int_0^T\langle F,u\rangle
 =\int_0^T\langle\mathfrak a,u\rangle.
\]
Consequently,
\[
 c_*\limsup_{\delta\to0}\|u_\delta-u\|_{L^2(0,T;V)}^2
 \leq\limsup_{\delta\to0}\int_0^T\bigl(\langle \mathcal Au_\delta,u_\delta\rangle
 -\langle \mathcal Au_\delta,u\rangle-\langle \mathcal Au,u_\delta-u\rangle\bigr)\dd t
 \leq\int_0^T\langle\mathfrak a,u\rangle-\int_0^T\langle\mathfrak a,u\rangle=0,
\]
since $\mathcal Au_\delta\weakto\mathfrak a$ and $u_\delta\weakto u$. Thus
$u_\delta\to u$ in $L^2(0,T;V)$, and the Lipschitz bound gives
$\mathcal Au_\delta\to \mathcal Au$, hence $\mathfrak a=\mathcal Au$.

Apply Lemma~\ref{lem:mass-energy} to this equation and use the definition
of $\mathcal A_\eps$. This gives \eqref{eq:micro-energy} for every $t\in[0,T]$.
If $u^1,u^2$ are two solutions corresponding to $F^1,F^2$ and the same
initial mass, put $d=u^1-u^2$; then $Jd(0)=0$. The same identity and strong monotonicity give
\[
 \frac12\|Jd(t)\|_H^2+c_*\int_0^t\|d\|_V^2\dd r
 \leq\int_0^t\langle F^1-F^2,d\rangle\dd r.
\]
Young's inequality yields the stability estimate
\begin{equation*}
 \|Jd(t)\|_H^2+c_*\int_0^t\|d\|_V^2\dd r
 \leq c_*^{-1}\int_0^t\|F^1-F^2\|_{V'}^2\dd r.
\end{equation*}
Taking $F^1=F^2$ proves uniqueness in $L^2(0,T;V)$.
For the bound, the energy identity, strong monotonicity and
$\mathcal A0=0$ give
\[
 \frac12\|Ju(t)\|_H^2+c_*\int_0^t\|u\|_V^2\dd r
 \leq\frac12\|u^{\rm in}\|_H^2+\int_0^t\langle F,u\rangle\dd r,
\]
and Young's inequality, with $\|u^{\rm in}\|_H^2\leq\|b\|_{L^\infty(Y)}
\|u^{\rm in}\|_{L^2(\Om)}^2$, bounds the first two terms in
\eqref{eq:micro-bound} by $CM_0^2$, with $C$ independent of $\eps$, $f$ and
$u^{\rm in}$. Finally,
$(\mathcal Bu)'=F-\mathcal Au$ and the Lipschitz bound on $\mathcal A$ bound the last term; the
constants are uniform in $\eps$.

For \eqref{eq:uniform-bulk}, Lemmas~\ref{lem:weighted-control} and
\ref{lem:extension-poincare} give
$\|\bar u_\eps\|_{L^2(0,T;W^{1,q}_0(\Om))}\leq C\|u_\eps\|_{L^2(0,T;V_\eps)}$,
and $w_i\leq\|w_i\|_{L^\infty(Y)}$ controls the fluxes by the weighted gradients.
The surface terms in \eqref{eq:uniform-mass} are bounded by
\eqref{eq:scaled-trace} and $|h(y,r)|\leq L_h|r|$. Let
$R_\eps:H^1_0(\Om)\to V_\eps$ be the restriction operator, which is
uniformly bounded since the weights are bounded. In $H^{-1}(\Om)$ we have
$m_\eps=R_\eps^*(\mathcal B_\eps u_\eps)$, and testing the equation with
$R_\eps\zeta$ gives, for a.e. $t$,
\begin{equation*}
 \langle m_\eps'(t),\zeta\rangle
 =\int_{\Omeps}f\zeta\dd x
   -\sum_{i=1}^N\int_{\Omeps}w_i^\eps D_i u_\eps D_i\zeta\dd x
 -\eps\int_{S_\eps}h(x/\eps,u_\eps)\zeta\dd\sigma_x
 \quad\forall\,\zeta\in H^1_0(\Om).
\end{equation*}
The unweighted form of \eqref{eq:scaled-trace} gives
$\eps\int_{S_\eps}|\zeta|^2\leq C\|\zeta\|_{H^1_0(\Om)}^2$, so all terms
are bounded and $m_\eps'$ is bounded in $L^2(0,T;H^{-1}(\Om))$.
Finally $\|m_\eps(t)\|_{L^2(\Om)}^2\leq\|b\|_{L^\infty(Y)}
\|J_\eps u_\eps(t)\|_{H_{b,\eps}}^2$.
\end{proof}

\begin{remark}\label{rem:zero-storage}
The initial condition prescribes only the mass, since $u_\eps(0)$ is not
defined in $V_\eps$. On the zero-storage component, the equation reduces to
the instantaneous constraints
$\langle \mathcal A_\eps u_\eps(t)-F_\eps(t),z\rangle=0$ for all $z\in\mathcal K_\eps$
and a.e. $t$; when $b\equiv0$, it is the elliptic equation
$\mathcal A_\eps u_\eps(t)=F_\eps(t)$ for a.e. $t$.
\end{remark}

\section{Strong compactness of the extended states}\label{sec:strong-compactness}

Theorem~\ref{thm:micro-solvability} controls the time derivative of the
mass $m_\eps$, not of the state. We first derive time-translation estimates
for the extended states, then combine them with Simon's compactness
criterion to obtain strong compactness of the states and masses, and finally
identify the weighted gradient limits.

\subsection{Detection and time translations}\label{sec:uniform}

To transfer the control of the mass to the extended states, we recover
$\bar u_\eps$ from the mass up to small errors; this is the first point of
the argument where the positivity of $\overline b$ becomes essential.

\begin{lemma}\label{lem:cell-detection}
Under $(H_1)$--$(H_3)$, let
\[
 q^*=\frac{Nq}{N-q}>2,\qquad
 \beta=\frac12-\frac1{q^*}=\frac1N-\frac1{2s}
       \geq\frac1{2N}>0.
\]
For all $v\in V_\eps$,
\begin{equation}
 \|P_\eps v\|_{L^2(\Om)}
 \leq C\bigl(\|\sqrt{b^\eps}\,v\|_{L^2(\Omeps)}
       +\eps^\beta\|v\|_{V_\eps}\bigr),
 \label{eq:cell-detection}
\end{equation}
with $C$ independent of $\eps$.
\end{lemma}

\begin{proof}
We first claim that
\begin{equation}
 \|z\|_{L^2(\Ystar)}
 \leq C\left(\|\sqrt b\,z\|_{L^2(\Ystar)}
            +\Bigl[\sum_{i=1}^N\int_{\Ystar}w_i|D_i z|^2\Bigr]^{1/2}\right)
 \qquad\forall\,z\in X_\omega(\Ystar).
 \label{eq:cell-weighted-mass}
\end{equation}
Otherwise, after normalization, there are $z_n$ with
$\|z_n\|_{L^2(\Ystar)}=1$ for which both terms on the right tend to
zero. Lemma~\ref{lem:weighted-control} bounds $z_n$ in
$W^{1,q}(\Ystar)$ and shows $Dz_n\to0$ in $L^q(\Ystar)^N$. Compactness of
$W^{1,q}(\Ystar)\hookrightarrow L^2(\Ystar)$ yields
$z_n\to c$ in $L^2(\Ystar)$, where $c$ is constant since $\Ystar$ is
connected. The first vanishing term then gives
$|c|^2\overline b=0$. Thus $c=0$, contrary to the normalization.

Let $E:X_\omega(\Ystar)\to W^{1,q}(Y)$ denote the cell
extension used in the proof of Lemma~\ref{lem:extension-poincare}. Its
continuity, the Sobolev embedding $W^{1,q}(Y)\hookrightarrow L^2(Y)$,
\eqref{eq:cell-weighted-mass}, and the weighted H\"older estimate give
\begin{equation}
 \|Ez\|_{L^2(Y)}
 \leq C\left(\|\sqrt b\,z\|_{L^2(\Ystar)}
        +\Bigl[\sum_{i=1}^N\int_{\Ystar}w_i|D_i z|^2\Bigr]^{1/2}\right)
 \qquad\forall\,z\in X_\omega(\Ystar).
 \label{eq:cell-extension-detection}
\end{equation}
For $k\in I_\eps$ and $z_k(y)=v(\eps(k+y))$, scaling
\eqref{eq:cell-extension-detection} and summing over complete cells gives
\begin{equation}
 \|P_\eps v\|_{L^2(\bigcup_{k\in I_\eps}\eps(k+Y))}
 \leq C\left(\|\sqrt{b^\eps}\,v\|_{L^2(\Omeps)}
      +\eps\Bigl[\sum_{i=1}^N\int_{\Omeps}w_i^\eps|D_i v|^2\Bigr]^{1/2}\right).
 \label{eq:cell-detection-interior}
\end{equation}
On the boundary strip, $|\Sigma_\eps|\leq C\eps$ by
Lemma~\ref{lem:periodic-averaging}; H\"older's inequality, the Sobolev
embedding and \eqref{eq:extension} give
\[
 \|P_\eps v\|_{L^2(\Sigma_\eps)}
 \leq |\Sigma_\eps|^\beta\|P_\eps v\|_{L^{q^*}(\Om)}
 \leq C\eps^\beta\|v\|_{V_\eps}.
\]
Since $\beta\leq1/2$ and $\eps<1$, the gradient term in
\eqref{eq:cell-detection-interior} is bounded by $\eps^\beta\|v\|_{V_\eps}$.
Combining the two estimates proves \eqref{eq:cell-detection}.
\end{proof}

\begin{remark}\label{rem:necessity}
The condition $\overline b>0$ is essential for the cell-detection
argument. Indeed, if $\overline b=0$, the nonnegativity of $b$ implies
$b=0$ a.e. on $\Ystar$. The microscopic storage is then confined to the
unperforated strip; for every $v\in V_\eps$, the extension bound,
$|\Sigma_\eps|\leq C\eps$, and $q^*>2$ give
\[
 \int_{\Omeps}b^\eps|v|^2\dd x
 =\int_{\Sigma_\eps}b^\eps|v|^2\dd x
 \leq C|\Sigma_\eps|^{1-2/q^*}
       \|P_\eps v\|_{L^{q^*}(\Om)}^2
 \leq C\eps^{2\beta}\|v\|_{V_\eps}^2.
\]
Thus uniformly energy-bounded states have vanishing stored mass in this
regime: by \eqref{eq:micro-bound}, the masses satisfy
$\|m_\eps\|_{L^2(Q)}\leq C\eps^\beta M_0\to0$, and the time derivative is
lost in the limit. This explains the role of positive cell-averaged storage
in the detection argument.
\end{remark}

Combined with the mass bound of Theorem~\ref{thm:micro-solvability}, the
detection inequality yields uniform control of the time translations.

\begin{proposition}\label{prop:time-translations}
Under $(H_1)$--$(H_5)$, for $0<\tau<T$,
\begin{align}
 &\int_0^{T-\tau}\int_{\Omeps}b^\eps
       |u_\eps(x,t+\tau)-u_\eps(x,t)|^2\dd x\dd t
       \leq C\tau M_0^2,\label{eq:weighted-translate}\\
 &\|\bar u_\eps(\cdot+\tau)-\bar u_\eps\|_
             {L^2(0,T-\tau;L^2(\Om))}
       \leq C(\sqrt\tau+\eps^\beta)M_0.
             \label{eq:state-translate}
\end{align}
In particular, the extended states have uniformly small time translations
along any sequence $\eps\to0$, in the sense that the $\limsup_{\eps\to0}$
of the left-hand side of \eqref{eq:state-translate} tends to zero as
$\tau\downarrow0$.
\end{proposition}

\begin{proof}
Since $\mathcal B_\eps$ is independent of time, Young's inequality for
convolution gives
\[
 \|\mathcal B_\eps[u_\eps(\cdot+\tau)-u_\eps]\|_
             {L^2(0,T-\tau;V_\eps')}
 \leq \tau\|(\mathcal B_\eps u_\eps)'\|_{L^2(0,T;V_\eps')}.
\]
Writing $d_\eps(t)=u_\eps(t+\tau)-u_\eps(t)$, we get
\begin{align*}
 \|\sqrt{b^\eps}\,d_\eps\|_{L^2(\Omeps\times(0,T-\tau))}^2
 &=\int_0^{T-\tau}\langle \mathcal B_\eps d_\eps,d_\eps\rangle\dd t\\
 &\leq \tau\|(\mathcal B_\eps u_\eps)'\|_{L^2(0,T;V_\eps')}
           \|d_\eps\|_{L^2(0,T-\tau;V_\eps)}\leq C\tau M_0^2,
\end{align*}
by \eqref{eq:micro-bound}. Finally apply
Lemma~\ref{lem:cell-detection} to $d_\eps(t)$, integrate in time,
and use \eqref{eq:weighted-translate} together with
$\|d_\eps\|_{L^2(0,T-\tau;V_\eps)}
 \leq2\|u_\eps\|_{L^2(0,T;V_\eps)}$.
Since $P_\eps$ is linear, this is \eqref{eq:state-translate}.
\end{proof}

\subsection{Strong convergence of the states and of the mass}

\begin{theorem}
\label{thm:strong-compactness}
Assume $(H_1)$--$(H_5)$, and let $\eps_n\to0$.
There exist a subsequence, not relabelled, and
\begin{equation}
 u_0\in L^2(0,T;W^{1,q}_0(\Om))\cap H^1(0,T;H^{-1}(\Om))
 \label{eq:compact-limit-space}
\end{equation}
such that
\begin{align}
 \bar u_{\eps_n}&\weakto u_0
       &&\text{in }L^2(0,T;W^{1,q}_0(\Om)),\label{eq:compact-weak}\\
 \bar u_{\eps_n}&\to u_0
       &&\text{strongly in }L^2(Q),\label{eq:compact-strong}\\
 \int_{\Omega_{\eps_n}\times(0,T)}|u_{\eps_n}-u_0|^2\dd x\dd t&\to0.
       &&\label{eq:compact-material}
\end{align}
Furthermore, with $\overline b>0$ as in $(H_3)$,
\begin{align}
 m_{\eps_n}&\weakto\overline b\,u_0
       &&\text{in }H^1(0,T;H^{-1}(\Om)),\label{eq:compact-mass-weak}\\
 m_{\eps_n}&\to\overline b\,u_0
       &&\text{in }C([0,T];H^{-1}(\Om)),\label{eq:compact-mass-strong}\\
 u_0(0)&=u^{\rm in}&&\text{in }H^{-1}(\Om).\label{eq:compact-initial}
\end{align}
In particular, $\overline b\,\partial_tu_0$ is the weak
$L^2(0,T;H^{-1}(\Om))$ limit of $m_{\eps_n}'$.
\end{theorem}

\begin{proof}
By \eqref{eq:uniform-bulk}, $(\bar u_{\eps_n})$ is bounded in
$L^2(0,T;W^{1,q}_0(\Om))$. Extract a subsequence such that
\eqref{eq:compact-weak} holds for some $u_0$ in that space.
To prove strong convergence, write
\[
 d_{n,\tau}(t)=\bar u_{\eps_n}(t+\tau)-\bar u_{\eps_n}(t)
 \quad(0<t<T-\tau).
\]
Proposition~\ref{prop:time-translations} gives
\begin{equation*}
 \|d_{n,\tau}\|_{L^2(0,T-\tau;L^2(\Om))}
 \leq C(\sqrt\tau+\eps_n^\beta)M_0.
\end{equation*}
The translations are therefore uniformly small. Indeed, given $\kappa>0$, choose $n_0$ so
large that $C\eps_n^\beta M_0<\kappa/3$
for $n\geq n_0$. Then choose $\tau_0>0$ sufficiently small that
$C\sqrt\tau M_0<\kappa/3$ for $\tau<\tau_0$.
For each of the finitely many indices $n<n_0$, translation continuity
of the fixed function $\bar u_{\eps_n}\in L^2(Q)$ gives the same bound
(after reducing $\tau_0$ if necessary). Thus
\begin{equation}
 \lim_{\tau\downarrow0}\sup_{n\geq1}
 \|d_{n,\tau}\|_{L^2(0,T-\tau;L^2(\Om))}=0.
 \label{eq:compact-uniform-translations}
\end{equation}
Because $q>2N/(N+2)$, the embedding
$W^{1,q}_0(\Om)\hookrightarrow L^2(\Om)$ is compact.
The time-translation compactness criterion of Simon
\cite{Simon1986}, applied to \eqref{eq:uniform-bulk} and
\eqref{eq:compact-uniform-translations}, yields relative compactness in
$L^2(Q)$. Its strong limit coincides with the weak limit in
\eqref{eq:compact-weak}, proving \eqref{eq:compact-strong}.
Since $P_{\eps_n}u_{\eps_n}=u_{\eps_n}$ on $\Omega_{\eps_n}$,
\[
 \int_{\Omega_{\eps_n}\times(0,T)}|u_{\eps_n}-u_0|^2
 \leq\|\bar u_{\eps_n}-u_0\|_{L^2(Q)}^2\to0,
\]
which gives \eqref{eq:compact-material}.

We next identify the limit of the extended mass.
By \eqref{eq:periodic-material-average} with $a=b$,
\begin{equation}
 \rho_\eps\weakto^*\overline b\quad\text{in }L^\infty(\Om).
 \label{eq:compact-rho-average}
\end{equation}

Decompose
\begin{equation*}
 m_{\eps_n}-\overline b\,u_0
 =\rho_{\eps_n}(\bar u_{\eps_n}-u_0)
   +(\rho_{\eps_n}-\overline b)u_0.
\end{equation*}
The first term converges to zero in $L^2(Q)$ by
\eqref{eq:compact-strong} and the uniform $L^\infty$ bound
on $\rho_\eps$. To treat the second, fix a time $t$ for which
$u_0(t)\in L^2(\Om)$. Equation~\eqref{eq:compact-rho-average}
shows that $(\rho_{\eps_n}-\overline b)u_0(t)$ tends weakly to zero
in $L^2(\Om)$. The compact embedding
$L^2(\Om)\hookrightarrow H^{-1}(\Om)$ then gives
\[
 \|(\rho_{\eps_n}-\overline b)u_0(t)\|_{H^{-1}(\Om)}
 \to0.
\]
The expression is bounded by $C\|u_0(t)\|_{L^2(\Om)}$, so dominated
convergence in time implies
\begin{equation}
 m_{\eps_n}\to\overline b\,u_0
 \quad\text{in }L^2(0,T;H^{-1}(\Om)).
 \label{eq:compact-mass-L2}
\end{equation}

The mass estimate \eqref{eq:uniform-mass} gives, after extraction,
$m_{\eps_n}\weakto m$ in $H^1(0,T;H^{-1}(\Om))$.
The strong limit in \eqref{eq:compact-mass-L2} identifies
$m=\overline b\,u_0$. Since $\overline b>0$,
$u_0\in H^1(0,T;H^{-1}(\Om))$, and
\eqref{eq:compact-mass-weak} follows. Furthermore, if
$v_n=m_{\eps_n}-\overline b\,u_0$, the Hilbert-space estimate
\[
 \|v_n\|_{C([0,T];H^{-1}(\Om))}^2
 \leq T^{-1}\|v_n\|_{L^2(0,T;H^{-1}(\Om))}^2
 +2\|v_n\|_{L^2(0,T;H^{-1}(\Om))}
      \|v_n'\|_{L^2(0,T;H^{-1}(\Om))}
\]
proves \eqref{eq:compact-mass-strong}, because $v_n'$ is bounded.
Finally, $m_{\eps_n}(0)=\rho_{\eps_n}u^{\rm in}$ by \eqref{eq:weak-micro} and
the restriction identity used in the proof of
Theorem~\ref{thm:micro-solvability}. By
\eqref{eq:periodic-material-average}, it converges weakly in $L^2(\Om)$ to
$\overline b\,u^{\rm in}$, hence strongly in $H^{-1}(\Om)$. Taking $t=0$ in
\eqref{eq:compact-mass-strong} gives $\overline b\,u_0(0)=\overline b\,u^{\rm in}$,
hence \eqref{eq:compact-initial}.
\end{proof}

\subsection{Two-scale compactness of the weighted gradients}

The extension estimates give a uniform bound in
$L^2(0,T;W^{1,q}_0(\Om))$ with $q<2$.
An affine coercivity inequality on the material cell allows us to recover
the $H^1$ regularity of the limit. We use
spatial two-scale convergence in the sense of
\cite{Allaire1992,Nguetseng1989}, with $(x,t)$ as the slow variables:
$z_\eps\overset{2}{\rightharpoonup}z$ in $L^p(Q)$ means
\[
 \int_Q z_\eps(x,t)\Theta(x,t,x/\eps)\dd x\dd t
 \to \int_Q\int_Y z(x,t,y)\Theta(x,t,y)\dd y\dd x\dd t
 \qquad\forall\,\Theta\in C_c^\infty\bigl(Q;C^\infty_\#(Y)\bigr).
\]
The boundedness of $z_\eps$ in $L^p(Q)$, with
$1<p<\infty$, ensures subsequential two-scale compactness, with
limits in $L^p(Q\times Y)$. There is no fast time
variable in the coefficients of \eqref{eq:micro}.

Let $\mathcal V_\omega$ be the completion of
\[
 \left\{\left.z\right|_{\Ystar}:z\in C^\infty_\#(Y),\
                  \int_{\Ystar}z\dd y=0\right\}
\]
for the norm
\begin{equation*}
 \|z\|_{\mathcal V_\omega}^2
    =\sum_{i=1}^N\int_{\Ystar}w_i|D_{y_i}z|^2\dd y.
\end{equation*}
The weighted H\"older inequality and Poincar\'e's inequality on the
connected cell show that this norm controls
$\|z\|_{W^{1,q}(\Ystar)}$ and $\|z\|_{L^2(\Ystar)}$; hence
$\mathcal V_\omega$ is a Hilbert space continuously embedded in
$W^{1,q}(\Ystar)$, consisting of functions with zero mean over $\Ystar$.
Since $\overline{\mathcal O}\subset Y$, the set $\Ystar$ contains a
neighbourhood of $\partial Y$ in $Y$, so the trace
$W^{1,q}(\Ystar)\to L^q(\partial Y)$ is continuous. For $j=1,\ldots,N$, the
condition that the traces of $z$ on the faces $\{y_j=0\}$ and $\{y_j=1\}$
coincide after translation by $e_j$ holds for restrictions of smooth periodic
functions and defines a closed linear subspace $W^{1,q}_{\rm per}(\Ystar)$ of
$W^{1,q}(\Ystar)$. This subspace is therefore weakly closed and contains
$\mathcal V_\omega$. In particular, a function of $W^{1,q}_{\rm per}(\Ystar)$
with constant gradient on $\Ystar$ has zero gradient.

\begin{lemma}\label{lem:periodic-gradient-representation}
Let $1<r<\infty$ and $K\in L^r_\#(Y)^N$. If
\begin{equation}
 \int_Y K(y)\cdot\Phi(y)\dd y=0
 \quad\forall\,\Phi\in C^\infty_\#(Y)^N
           \text{ with }\Div_y\Phi=0,
 \label{eq:periodic-gradient-orthogonality}
\end{equation}
then there is a unique $v\in W^{1,r}_\#(Y)$ of zero cell mean such that
$K=D_yv$. Moreover,
\begin{equation}
 \|v\|_{W^{1,r}(Y)}\leq C_r\|K\|_{L^r(Y)^N}.
 \label{eq:periodic-gradient-bound}
\end{equation}
The map $K\mapsto v$ on the closed subspace defined by
\eqref{eq:periodic-gradient-orthogonality} is linear and bounded.
\end{lemma}

\begin{proof}
For $r=2$ the statement is the Weyl decomposition
\cite[Section~1.1]{JikovKozlovOleinik1994}; the following argument covers
$1<r<\infty$.

Constant vector fields are admissible in
\eqref{eq:periodic-gradient-orthogonality}; taking $\Phi=e_i$ for
$i=1,\ldots,N$ gives $\int_YK\dd y=0$. Next, fix $i\neq j$ and
$\zeta\in C^\infty_\#(Y)$, and let $\Phi$ have components
$\Phi_i=D_{y_j}\zeta$, $\Phi_j=-D_{y_i}\zeta$ and $\Phi_l=0$ for
$l\notin\{i,j\}$. Then $\Div_y\Phi=D_{y_i}D_{y_j}\zeta-D_{y_j}D_{y_i}\zeta=0$,
and \eqref{eq:periodic-gradient-orthogonality} gives
\[
 0=\int_Y\bigl(K_iD_{y_j}\zeta-K_jD_{y_i}\zeta\bigr)\dd y
 =\langle D_{y_i}K_j-D_{y_j}K_i,\zeta\rangle_{\mathscr D'_\#(Y),C^\infty_\#(Y)}.
\]
Hence $D_{y_j}K_i=D_{y_i}K_j$ in $\mathscr D'_\#(Y)$ for all $i,j$.

Let $(\eta_\delta)_{\delta>0}$ be a standard mollifier and
$K_\delta=K*\eta_\delta$. Since $K$ is $Y$-periodic, $K_\delta$ is smooth
and $Y$-periodic; convolution commutes with differentiation and preserves
the cell mean, so $\int_YK_\delta\dd y=0$ and $D_{y_j}(K_\delta)_i=D_{y_i}(K_\delta)_j$.
Moreover $K_\delta\to K$ in $L^r(Y)^N$, because $r<\infty$.

We next construct a potential of $K_\delta$ by Fourier series. Write
\[
 K_\delta(y)=\sum_{k\in\mathbb Z^N\setminus\{0\}}\widehat K_\delta(k)
 e^{2\pi\mathrm ik\cdot y},
\]
where the term $k=0$ is absent because $K_\delta$ has zero mean, and the
coefficients decay faster than any power of $|k|$ because $K_\delta$ is
smooth. The curl-free relations become
$k_j\widehat K_{\delta,i}(k)=k_i\widehat K_{\delta,j}(k)$ for all $i,j$ and
$k\neq0$. Multiplying by $k_j$ and summing over $j$ gives
$|k|^2\widehat K_{\delta,i}(k)=k_i\,\bigl(k\cdot\widehat K_\delta(k)\bigr)$, that
is, $\widehat K_\delta(k)=c_k k$ with
$c_k=|k|^{-2}k\cdot\widehat K_\delta(k)$. Set
\[
 v_\delta(y)=\sum_{k\in\mathbb Z^N\setminus\{0\}}
 \frac{c_k}{2\pi\mathrm i}\,e^{2\pi\mathrm ik\cdot y}.
\]
The coefficients decay rapidly, so $v_\delta\in C^\infty_\#(Y)$; it has zero
cell mean, and termwise differentiation gives $D_yv_\delta=K_\delta$ on $Y$.
Since $K_\delta$ is real, $\widehat K_\delta(-k)=\overline{\widehat K_\delta(k)}$,
hence $c_{-k}=-\overline{c_k}$; the coefficients
$c_k/(2\pi\mathrm i)$ of $v_\delta$ therefore satisfy the same conjugation
relation, and $v_\delta$ is real-valued.

The Poincar\'e--Wirtinger inequality
$\|z\|_{L^r(Y)}\leq C_r'\|D_yz\|_{L^r(Y)^N}$ for mean-zero
$z\in W^{1,r}_\#(Y)$, applied to $z=v_\delta-v_{\delta'}$, gives
\[
 \|v_\delta-v_{\delta'}\|_{W^{1,r}(Y)}
 \leq(1+C_r')\|K_\delta-K_{\delta'}\|_{L^r(Y)^N}.
\]
Thus $(v_\delta)$ is a Cauchy sequence as $\delta\to0$. Its limit $v$ belongs to
$W^{1,r}_\#(Y)$, since periodicity is preserved under convergence in
$W^{1,r}_{\rm loc}(\R^N)$; it has zero cell mean and satisfies $D_yv=K$.
Passing to the limit in $\|v_\delta\|_{W^{1,r}(Y)}\leq(1+C_r')
\|K_\delta\|_{L^r(Y)^N}$ gives \eqref{eq:periodic-gradient-bound} with
$C_r=1+C_r'$. If $v$ and $\tilde v$ are two such potentials, then
$D_y(v-\tilde v)=0$ on the connected torus, so $v-\tilde v$ is constant; the
zero-mean normalization makes it vanish. This proves uniqueness.

By uniqueness, the map $\Upsilon:K\mapsto v$ is linear on the closed subspace of
$L^r_\#(Y)^N$ defined by \eqref{eq:periodic-gradient-orthogonality}, and
\eqref{eq:periodic-gradient-bound} shows that it is bounded into
$W^{1,r}_\#(Y)$. Consequently, if $K\in L^r(Q;L^r_\#(Y)^N)$ satisfies
\eqref{eq:periodic-gradient-orthogonality} for a.e. $(x,t)\in Q$, then
$(x,t)\mapsto\Upsilon K(x,t,\cdot)$ is strongly measurable, as the composition of a
bounded linear operator with a strongly measurable function, and
$\|\Upsilon K\|_{L^r(Q;W^{1,r}_\#(Y))}\leq C_r\|K\|_{L^r(Q\times Y)^N}$.
\end{proof}

\begin{lemma}\label{lem:twoscale-affine}
Under $(H_1)$--$(H_2)$, there exists $c_\omega>0$ such that
\begin{equation}
 c_\omega\bigl(|\xi|^2+\|z\|_{\mathcal V_\omega}^2\bigr)
 \leq\sum_{i=1}^N\int_{\Ystar}w_i
                  |\xi_i+D_{y_i}z|^2\dd y
 \quad\forall\,\xi\in\R^N,\ \forall\,z\in\mathcal V_\omega.
 \label{eq:twoscale-affine}
\end{equation}
In particular, the map
\begin{equation}
 \mathcal T_\omega(\xi,z)
 =\bigl(\sqrt{w_i}(\xi_i+D_{y_i}z)\bigr)_{i=1}^N,
 \quad\mathcal T_\omega:\R^N\times\mathcal V_\omega
             \to L^2(\Ystar)^N,
 \label{eq:twoscale-T}
\end{equation}
is injective and has closed range $\mathcal M_\omega$.
\end{lemma}

\begin{proof}
We first prove that
\begin{equation}
 |\xi|\leq C\|\mathcal T_\omega(\xi,z)\|_{L^2(\Ystar)^N}
 \qquad\forall\,(\xi,z)\in\R^N\times\mathcal V_\omega.
 \label{eq:twoscale-xi-control}
\end{equation}
Otherwise one can choose $|\xi_n|=1$ and $z_n\in\mathcal V_\omega$
with $\|\mathcal T_\omega(\xi_n,z_n)\|_{L^2(\Ystar)^N}\to0$.
Since $w_i\in L^\infty(Y)$,
\[
 \sum_{i=1}^N\|\sqrt{w_i}D_i z_n\|_{L^2(\Ystar)}^2
 \leq 2\|\mathcal T_\omega(\xi_n,z_n)\|_{L^2(\Ystar)^N}^2
       +2\sum_{i=1}^N\int_{\Ystar}w_i|\xi_{n,i}|^2\dd y\leq C.
\]
The weighted H\"older estimate and zero-mean Poincar\'e inequality
bound $(z_n)$ in $W^{1,q}(\Ystar)$. After extraction,
$\xi_n\to\xi$ with $|\xi|=1$, and
$z_n\rightharpoonup z$ in $W^{1,q}(\Ystar)$; since $W^{1,q}_{\rm per}(\Ystar)$
is weakly closed, $z\in W^{1,q}_{\rm per}(\Ystar)$.
Furthermore, since $q/(2-q)=s$, H\"older's inequality gives
\[
 \|\xi_n+D_yz_n\|_{L^q(\Ystar)}
 \leq C\|\mathcal T_\omega(\xi_n,z_n)\|_{L^2(\Ystar)^N}
 \to0.
\]
Thus $D_yz=-\xi$ in $\Ystar$. Connectedness implies
$z(y)=-\xi\cdot y+c$. Since the traces of $z$ coincide on
opposite faces of $Y$, all components of $\xi$ vanish, a
contradiction. This proves \eqref{eq:twoscale-xi-control}.
Finally,
\[
 \|z\|_{\mathcal V_\omega}^2
 \leq2\|\mathcal T_\omega(\xi,z)\|_{L^2(\Ystar)^N}^2
     +2\sum_{i=1}^N\int_{\Ystar}w_i|\xi_i|^2\dd y
 \leq C\|\mathcal T_\omega(\xi,z)\|_{L^2(\Ystar)^N}^2.
\]
The two bounds give \eqref{eq:twoscale-affine}. The map
\eqref{eq:twoscale-T} is bounded; its inverse on its range is
bounded by \eqref{eq:twoscale-affine}. Completeness of
$\R^N\times\mathcal V_\omega$ implies that the range is closed.
\end{proof}

The closed range allows us to identify the two-scale limit of the weighted
gradients and to recover the $H^1$ regularity of $u_0$.

\begin{theorem}
\label{thm:twoscale-weighted}
Assume $(H_1)$--$(H_5)$ and let $\eps_n\to0$. After a further
extraction from Theorem~\ref{thm:strong-compactness}, there exists a
unique corrector
\[
 u_1\in L^2(Q;\mathcal V_\omega)
\]
associated with the chosen two-scale limit such that
\begin{equation}
 u_0\in L^2(0,T;H^1_0(\Om))\cap H^1(0,T;H^{-1}(\Om))
 \label{eq:twoscale-u0-H1}
\end{equation}
and, for $i=1,\ldots,N$,
\begin{align}
 \chi_{\Omega_{\eps_n}}D_i u_{\eps_n}
 &\overset{2}{\rightharpoonup}
 \chi_{\Ystar}(D_{x_i}u_0+D_{y_i}u_1)
       &&\text{in }L^q(Q),\label{eq:twoscale-gradient-q}\\
 \chi_{\Omega_{\eps_n}}\sqrt{w_i^{\eps_n}}D_i u_{\eps_n}
 &\overset{2}{\rightharpoonup}
 \chi_{\Ystar}\sqrt{w_i}\,(D_{x_i}u_0+D_{y_i}u_1)
       &&\text{in }L^2(Q),\label{eq:twoscale-gradient-weighted}\\
 \chi_{\Omega_{\eps_n}}w_i^{\eps_n}D_i u_{\eps_n}
 &\overset{2}{\rightharpoonup}
 \chi_{\Ystar}w_i(D_{x_i}u_0+D_{y_i}u_1)
       &&\text{in }L^2(Q).\label{eq:twoscale-flux}
\end{align}
Here the limits are extended by zero in $\mathcal O$ and understood on
$Q\times Y$. In particular, the zero-extended physical fluxes obey
\begin{equation}
 \chi_{\Omega_{\eps_n}}w_i^{\eps_n}D_i u_{\eps_n}
 \rightharpoonup \int_{\Ystar}w_i
 (D_{x_i}u_0+D_{y_i}u_1)\dd y
 \quad\text{weakly in }L^2(Q).
 \label{eq:twoscale-flux-average}
\end{equation}
In addition,
\begin{align}
 &\|D_xu_0\|_{L^2(Q)^N}^2
 +\|u_1\|_{L^2(Q;\mathcal V_\omega)}^2
       \leq CM_0^2,\label{eq:twoscale-limit-bound}\\
 &\int_Q\int_{\Ystar}\sum_{i=1}^{N}w_i
       |D_{x_i}u_0+D_{y_i}u_1|^2\dd y\dd x\dd t
 \leq\liminf_{n\to\infty}\sum_{i=1}^N\int_{\Omega_{\eps_n}\times(0,T)}
           w_i^{\eps_n}|D_i u_{\eps_n}|^2\dd x\dd t.
 \label{eq:twoscale-liminf}
\end{align}
\end{theorem}

\begin{proof}
By \eqref{eq:uniform-bulk}, $D\bar u_{\eps_n}$ is bounded in
$L^q(Q)^N$. The two-scale compactness theorem \cite{Allaire1992,Nguetseng1989} provides
$\Lambda\in L^q(Q\times Y)^N$ as its two-scale limit. Since the two-scale
limit of $\bar u_{\eps_n}$ is $u_0$ by the strong convergence
\eqref{eq:compact-strong}, integration by parts against any
$\psi\in C_c^\infty\bigl(Q;C^\infty_\#(Y)^N\bigr)$ with $\Div_y\psi=0$ gives
\[
 \int_Q\int_Y\Lambda\cdot\psi\dd y\dd x\dd t
 =\int_Q\int_YD_xu_0\cdot\psi\dd y\dd x\dd t.
\]
To localize this orthogonality, take tests of the form
$\psi(x,t,y)=\phi(x,t)\Phi(y)$ with
$\phi\in C_c^\infty(Q)$ and $\Phi\in C^\infty_\#(Y)^N$
solenoidal. The distributional fundamental lemma and a countable
dense family of such $\Phi$ show that, for a.e. $(x,t)$,
$K(x,t,\cdot):=\Lambda(x,t,\cdot)-D_xu_0(x,t)$ satisfies
\eqref{eq:periodic-gradient-orthogonality}. Apply
Lemma~\ref{lem:periodic-gradient-representation} with $r=q$ to this
field. Its bounded linear potential map is Bochner measurable, and we
obtain a unique mean-zero $\widetilde u_1\in L^q(Q;W^{1,q}_\#(Y))$ with
$\|\widetilde u_1\|_{L^q(Q;W^{1,q}_\#(Y))}\leq C\|K\|_{L^q(Q\times Y)^N}$ such
that
\begin{equation*}
 D_i\bar u_{\eps_n}\overset{2}{\rightharpoonup}
 D_{x_i}u_0+D_{y_i}\widetilde u_1
 \quad\text{in }L^q(Q),\qquad i=1,\ldots,N.
\end{equation*}

Since no holes are removed in the strip,
\[
 \chi_{\Omeps}D_i u_\eps-\chi_{\Ystar}(x/\eps)D_i\bar u_\eps
 =\mathbf 1_{\Sigma_\eps}\bigl(1-\chi_{\Ystar}(x/\eps)\bigr)D_i\bar u_\eps,
\]
which is bounded in $L^q(Q)$ and has zero two-scale limit by
\eqref{eq:strip-error}. The multiplier statement \eqref{eq:periodic-multiplier} of
Lemma~\ref{lem:periodic-averaging}, with $p=q$ and $a=\chi_{\Ystar}$,
therefore identifies the two-scale limit of the material gradients as
\begin{equation}
 \chi_{\Omeps}D_i u_\eps\overset{2}{\rightharpoonup}
 \chi_{\Ystar}(D_{x_i}u_0+D_{y_i}\widetilde u_1)
 \quad\text{in }L^q(Q).
 \label{eq:twoscale-gradient-preliminary}
\end{equation}
By \eqref{eq:micro-bound}, the weighted gradients
$g_{\eps,i}:=\chi_{\Omeps}\sqrt{w_i^\eps}D_i u_\eps$
are bounded in $L^2(Q)$; extracting again, $g_{\eps,i}$ converges weakly
two-scale to some $G_i$ in $L^2(Q)$. Since $\bar u_\eps=u_\eps$ on the strip,
\[
 g_{\eps,i}-D_i\bar u_\eps(\chi_{\Ystar}\sqrt{w_i})(x/\eps)
 =\mathbf 1_{\Sigma_\eps}\bigl(1-\chi_{\Ystar}(x/\eps)\bigr)g_{\eps,i},
\]
which is bounded in $L^2(Q)$ and has zero two-scale limit by
\eqref{eq:strip-error}. Applying \eqref{eq:periodic-multiplier}, with $p=q$ and
$a=\chi_{\Ystar}\sqrt{w_i}\in L^\infty_\#(Y)$, gives
\begin{equation}
 G_i=\chi_{\Ystar}\sqrt{w_i}\,
       (D_{x_i}u_0+D_{y_i}\widetilde u_1)
 \quad\text{a.e. on }Q\times Y.
 \label{eq:twoscale-G-preliminary}
\end{equation}
In particular, $G_i=0$ on $Q\times\mathcal O$.

We next exploit the closed range $\mathcal M_\omega$ in
Lemma~\ref{lem:twoscale-affine}. Let
$\Xi\in\mathcal M_\omega^\perp\subset L^2(\Ystar)^N$, and set
$\Psi_i=\sqrt{w_i}\,\Xi_i$ on $\Ystar$. Extend both $\Xi$ and $\Psi$
by zero into $\mathcal O$ and periodically to $\R^N$. Orthogonality to
$\mathcal T_\omega(\xi,0)$ and
$\mathcal T_\omega(0,z)$ implies, respectively,
\begin{equation*}
 \int_Y\Psi(y)\dd y=0,
 \qquad \int_Y\Psi\cdot D_y\zeta\dd y=0
           \quad\forall\,\zeta\in C^\infty_\#(Y).
\end{equation*}
Indeed, the restriction of any $\zeta\in C^\infty_\#(Y)$,
after subtracting its mean over $\Ystar$, lies in
$\mathcal V_\omega$, and $\Psi=0$ in $\mathcal O$. Consequently
\[
 \langle\Div_y\Psi,\zeta\rangle_{\mathscr D'_\#(Y),
                                    C^\infty_\#(Y)}
 =-\int_{\Ystar}\Psi\cdot D_y\zeta\dd y=0.
\]
This equality holds on the \emph{whole torus}, including for tests
meeting $S$: the distributional divergence of the zero extension
vanishes across the hole boundary as well as inside the material. In particular,
$\Psi\in L^2_\#(Y)^N$ satisfies
$\Div_y\Psi=0$ in $\mathscr D'_\#(Y)$.
Periodic rescaling, first against smooth tests and then by density,
gives $\Div_x[\Psi(x/\eps)]=0$ in $\mathscr D'(\R^N)$.

For every $\phi\in C_c^\infty(Q)$, the distributional Green formula gives
\begin{equation}
 \int_Q\phi\sum_{i=1}^N g_{\eps,i}\Xi_i(x/\eps)\dd x\dd t
 =-\int_Q\bar u_\eps \Psi(x/\eps)\cdot D_x\phi\dd x\dd t.
 \label{eq:twoscale-Green}
\end{equation}
For smooth states, extend the state smoothly across the complete
holes and apply the distributional divergence identity to its product
with $\phi$. The integrand $\Psi(x/\eps)$ vanishes on the
complete holes, so their contributions vanish identically. On the
unperforated boundary strip we retain the original state and its
gradient. The compact support of $\phi$ removes the
exterior boundary contribution, which gives \eqref{eq:twoscale-Green}.
For a general state, fix $\eps$. The defining dense class of
$H^1_{\omega^\eps}(\Omeps)$ consists of restrictions of functions
smooth in a neighbourhood of $\overline{\Omeps}$. Finite sums of
simple time functions with values in that class are dense in the
Bochner space $L^2(0,T;H^1_{\omega^\eps}(\Omeps))$; each time
coefficient can in turn be approximated by a smooth compactly
supported time function. We may thus select $v_n$, smooth in space
near $\overline{\Omeps}$ and smooth in time, converging to
$u_\eps$ in this Bochner norm. As the number of holes is finite for
fixed $\eps$, each spatial approximant extends smoothly through the
holes; denote one such extension by $\widetilde v_n$. The identity
$\Div_x\Psi(x/\eps)=0$ can then be tested with
$\phi\widetilde v_n$ by ordinary distributional integration by
parts. Both integrands vanish within complete holes because
$\Psi(x/\eps)=0$ there, so the identity depends only on $v_n|_{\Omeps}$;
since $\phi$ has compact support, the exterior trace of $v_n$ does not enter.
Its left-hand side converges,
since
\[
 \left|\int_{\Qeps}\phi\sum_{i=1}^N
   \sqrt{w_i^\eps}D_i(v_n-u_\eps)\,\Xi_i(x/\eps)\right|
 \leq\|\phi\|_{L^\infty(Q)}
 \left(\sum_{i=1}^N\|\sqrt{w_i^\eps}D_i(v_n-u_\eps)\|_{L^2(\Qeps)}^2\right)^{1/2}
 \|\Xi(x/\eps)\|_{L^2(\Om)^N}.
\]
The right-hand side converges directly on $\Qeps$ by the
$L^2(\Qeps)$ convergence of $v_n$ and
$\Psi(x/\eps)\in L^2(\Om)^N$. Finally,
$\Psi(x/\eps)=0$ in $\Om\setminus\Omeps$ and
$P_\eps u_\eps=u_\eps$ on $\Omeps$, so the resulting identity
is exactly \eqref{eq:twoscale-Green}.

To pass to the limit, apply \eqref{eq:periodic-multiplier} with $p=2$,
$a=\Xi_i\in L^2_\#(Y)$ and $\Theta=\phi$ to the bounded sequence
$g_{\eps,i}$; the left-hand side tends to
$\int_Q\phi\int_{\Ystar}G\cdot\Xi$.
On the right, $\bar u_\eps\to u_0$ strongly in $L^2(Q)$ and
$\Psi(x/\eps)$ is uniformly bounded in $L^2(\Om)^N$ by
\eqref{eq:periodic-cell-average}; replacing $\bar u_\eps$ by $u_0$
therefore incurs an $o(1)$ error. Finally, approximate the $L^2(Q)$
function $u_0D_x\phi$ by smooth compactly supported functions,
apply the periodic mean-value formula to each approximation, and use
the same uniform $L^2(\Om)^N$ bound to remove it. The limit is zero since
$\int_Y\Psi=0$. Hence \eqref{eq:twoscale-Green} gives
\[
 \int_Q\phi(x,t)\int_{\Ystar}G(x,t,y)\cdot\Xi(y)
       \dd y\dd x\dd t=0.
\]
A countable dense subset of $\mathcal M_\omega^\perp$ shows that
$G(x,t,\cdot)\in\mathcal M_\omega$ for a.e. $(x,t)\in Q$.
By the bounded inverse of $\mathcal T_\omega$ on its range,
there are unique measurable
$\xi\in L^2(Q)^N$ and $z\in L^2(Q;\mathcal V_\omega)$
satisfying
\begin{align}
 G_i&=\chi_{\Ystar}\sqrt{w_i}\,(\xi_i+D_{y_i}z)\quad\text{a.e. on }Q\times Y,
                                      \label{eq:twoscale-affine-rep}\\
 \|\xi\|_{L^2(Q)}^2+\|z\|_{L^2(Q;\mathcal V_\omega)}^2
 &\leq C\|G\|_{L^2(Q\times Y)^N}^2.
                                      \label{eq:twoscale-affine-bound}
\end{align}

Compare \eqref{eq:twoscale-affine-rep} and
\eqref{eq:twoscale-G-preliminary}. Since $w_i>0$ a.e. on $\Ystar$,
\[
 D_y(z-\widetilde u_1)=D_xu_0-\xi
 \quad\text{a.e. on }Q\times\Ystar.
\]
For a.e. $(x,t)\in Q$, both $z(x,t,\cdot)\in\mathcal V_\omega$ and the
restriction of $\widetilde u_1(x,t,\cdot)\in W^{1,q}_\#(Y)$ to $\Ystar$ belong
to $W^{1,q}_{\rm per}(\Ystar)$. Their difference has a constant gradient on
the connected cell, which therefore vanishes. Hence $\xi=D_xu_0$ a.e. in $Q$. Choosing the
mean-zero representative of $\widetilde u_1|_{\Ystar}$ now gives
$z=\widetilde u_1|_{\Ystar}$. Put $u_1=z$.
Estimate \eqref{eq:twoscale-affine-bound} proves
$D_xu_0\in L^2(Q)^N$ and the bound for $u_1$.
Since the trace of $u_0$ on $\partial\Om$ vanishes by
\eqref{eq:compact-limit-space}, this establishes
\eqref{eq:twoscale-u0-H1} and \eqref{eq:twoscale-limit-bound}.
Equations \eqref{eq:twoscale-gradient-preliminary} and
\eqref{eq:twoscale-affine-rep} give
\eqref{eq:twoscale-gradient-q}--\eqref{eq:twoscale-gradient-weighted}.
Since $\chi_{\Omeps}w_i^\eps D_iu_\eps=g_{\eps,i}\sqrt{w_i}(x/\eps)$ is
bounded in $L^2(Q)$, the multiplier statement of
Lemma~\ref{lem:periodic-averaging}, with $p=r=2$ and $a=\sqrt{w_i}$, gives
\eqref{eq:twoscale-flux}.
Testing \eqref{eq:twoscale-flux} with functions independent of $y$
proves \eqref{eq:twoscale-flux-average}.
Finally, weak lower semicontinuity of the $L^2$ norm under
two-scale convergence proves \eqref{eq:twoscale-liminf}.
\end{proof}

\section{The cell problem and the homogenized equation}\label{sec:homogenization}

The effective diffusion is given by a stationary weighted Neumann cell
problem, which we solve first.

\subsection{The weighted cell problem and effective diffusion}

For each $\xi\in\R^N$, seek $\pi_\xi\in\mathcal V_\omega$ satisfying
\begin{equation}
 \sum_{i=1}^N\int_{\Ystar}w_i
      (\xi_i+D_{y_i}\pi_\xi)D_{y_i}z\dd y=0
       \qquad\forall\,z\in\mathcal V_\omega.
 \label{eq:hom-cell}
\end{equation}
The mean-zero normalization in the definition of
$\mathcal V_\omega$ fixes the additive constant. The equation is the
weak form of
\begin{equation}
 \left\{
 \begin{aligned}
 -\Div_y\bigl(A(\xi+D_y\pi_\xi)\bigr)&=0
                &&\text{in }\Ystar,\\
 A(\xi+D_y\pi_\xi)\cdot\nu_{\Ystar}&=0
                &&\text{on }S,
 \end{aligned}\right.
 \label{eq:hom-cell-formal}
\end{equation}
with periodic transmission across opposite cell faces, where
$\nu_{\Ystar}$ is the unit normal exterior to $\Ystar$. The conormal
condition in \eqref{eq:hom-cell-formal} is understood in the weak sense of
\eqref{eq:hom-cell}.

\begin{proposition}
\label{prop:hom-cell}
Under $(H_1)$--$(H_2)$, problem~\eqref{eq:hom-cell} has exactly one
solution for each $\xi\in\R^N$. The map $\xi\mapsto\pi_\xi$ is linear,
and
\begin{equation}
 \|\pi_\xi\|_{\mathcal V_\omega}\leq C|\xi|.
 \label{eq:hom-cell-bound}
\end{equation}
Let $\pi^j=\pi_{e_j}$ and define
\begin{equation}
 A^{\rm hom}_{ij}
   =\int_{\Ystar}w_i
      \bigl(\delta_{ij}+D_{y_i}\pi^j\bigr)\dd y,
          \qquad 1\leq i,j\leq N.
 \label{eq:hom-matrix}
\end{equation}
Then $A^{\rm hom}$ is a constant, symmetric, positive-definite matrix:
with $c_\omega$ from \eqref{eq:twoscale-affine},
\begin{equation}
 c_\omega|\xi|^2\leq A^{\rm hom}\xi\cdot\xi
 =\min_{z\in\mathcal V_\omega}\sum_{i=1}^N\int_{\Ystar}w_i|\xi_i+D_{y_i}z|^2\dd y
 \leq\sum_{i=1}^N\Bigl(\int_{\Ystar}w_i\dd y\Bigr)\xi_i^2
           \qquad\forall\,\xi\in\R^N.
 \label{eq:hom-ellipticity}
\end{equation}
More precisely,
\begin{equation}
 A^{\rm hom}\xi\cdot\eta
   =\sum_{i=1}^N\int_{\Ystar}w_i
   (\xi_i+D_{y_i}\pi_\xi)
   (\eta_i+D_{y_i}\pi_\eta)\dd y
        \qquad\forall\,\xi,\eta\in\R^N.
 \label{eq:hom-energy-identity}
\end{equation}
\end{proposition}

\begin{proof}
The form
\[
 a_\omega(z,v)=\sum_{i=1}^N\int_{\Ystar}
                           w_iD_{y_i}zD_{y_i}v\dd y,\qquad z,v\in\mathcal V_\omega,
\]
is an inner product defining the Hilbert norm of
$\mathcal V_\omega$. Moreover,
\[
 \left|\sum_{i=1}^N\int_{\Ystar}w_i\xi_iD_{y_i}z\right|
 \leq\left(\sum_{i=1}^N\int_{\Ystar}w_i|\xi_i|^2\right)^{1/2}
             \|z\|_{\mathcal V_\omega}
 \leq C|\xi|\|z\|_{\mathcal V_\omega}.
\]
Lax--Milgram therefore gives a unique $\pi_\xi$ solving
\eqref{eq:hom-cell}. Testing with $z=\pi_\xi$ proves
\eqref{eq:hom-cell-bound}; linearity follows from uniqueness.

Write $Z_\xi=A(\xi+D_y\pi_\xi)$. By definition,
$A^{\rm hom}\xi=\int_{\Ystar}Z_\xi\dd y$. Using
$z=\pi_\eta$ in \eqref{eq:hom-cell} gives
\begin{align*}
 A^{\rm hom}\xi\cdot\eta
 &=\sum_{i=1}^N\int_{\Ystar}w_i
       (\xi_i+D_{y_i}\pi_\xi)\eta_i\dd y\\
 &=\sum_{i=1}^N\int_{\Ystar}w_i
       (\xi_i+D_{y_i}\pi_\xi)
       (\eta_i+D_{y_i}\pi_\eta)\dd y.
\end{align*}
This is \eqref{eq:hom-energy-identity} and proves symmetry.
For $\eta=\xi$, the affine cell inequality
\eqref{eq:twoscale-affine} gives
\[
 A^{\rm hom}\xi\cdot\xi
  =\sum_{i=1}^N\int_{\Ystar}w_i|\xi_i+D_{y_i}\pi_\xi|^2\dd y
  \geq c_\omega|\xi|^2.
\]
Finally, \eqref{eq:hom-cell} is the Euler--Lagrange equation of the strictly
convex functional $z\mapsto\sum_{i=1}^N\int_{\Ystar}w_i|\xi_i+D_{y_i}z|^2\dd y$
on $\mathcal V_\omega$, whose minimizer is therefore $\pi_\xi$; this gives the
minimum characterization, and the choice $z=0$ gives the upper bound in
\eqref{eq:hom-ellipticity}.
\end{proof}

\subsection{Identification of the microscopic corrector}

We show that the corrector of Theorem~\ref{thm:twoscale-weighted} solves the
cell problem \eqref{eq:hom-cell}, which determines it uniquely.

\begin{proposition}
\label{prop:hom-corrector}
Under $(H_1)$--$(H_5)$, let $(u_0,u_1)$ be a subsequential limit
from Theorems~\ref{thm:strong-compactness} and
\ref{thm:twoscale-weighted}. Then, for a.e. $(x,t)\in Q$,
\begin{equation}
 u_1(x,t,\cdot)=\pi_{D_xu_0(x,t)}
             =\sum_{j=1}^N\pi^j(\cdot)D_{x_j}u_0(x,t)
             \quad\text{in }\mathcal V_\omega.
 \label{eq:hom-corrector-identification}
\end{equation}
Consequently, after extraction,
\begin{equation}
 \chi_{\Omeps}A^\eps Du_\eps
       \rightharpoonup A^{\rm hom}D_xu_0
                   \quad\text{weakly in }L^2(Q)^N.
 \label{eq:hom-flux-convergence}
\end{equation}
\end{proposition}

\begin{proof}
For $\varphi\in C_c^\infty(Q)$ and $\psi\in C^\infty_\#(Y)$, take in
\eqref{eq:weak-micro} the time-dependent test
\[
 v_\eps(x,t)=\eps\varphi(x,t)\psi(x/\eps)|_{\Omeps}.
\]
It belongs to $L^2(0,T;V_\eps)$ and is smooth in time.
Integrating in time and using compact support of $\varphi$ gives
\begin{align}
 &-\eps\int_{\Qeps}b^\eps u_\eps
                  \partial_t\varphi\,\psi(x/\eps)\dd x\dd t
  +\int_{\Qeps}A^\eps Du_\eps\cdot
     \bigl(\varphi D_y\psi(x/\eps)
              +\eps\psi(x/\eps)D_x\varphi\bigr)\dd x\dd t
       \nonumber\\
 &\quad+\eps^2\int_0^T\int_{S_\eps}
            h(x/\eps,u_\eps)\varphi\psi(x/\eps)
                           \dd\sigma_x\dd t
   =\eps\int_{\Qeps}f\varphi\psi(x/\eps)\dd x\dd t.
 \label{eq:hom-oscillating-test}
\end{align}
The first term and the right-hand side are $O(\eps)$ by the
bounds on $u_\eps$ in $L^2(\Qeps)$ and on $f$ in $L^2(Q)$ and boundedness of $b$.
The term containing $\eps D_x\varphi$ is $O(\eps)$ by
\eqref{eq:uniform-bulk}. To estimate the surface term, use
Cauchy--Schwarz in the measure $\eps\dd\sigma_x\dd t$:
\begin{align*}
 &\left|\eps^2\int_0^T\int_{S_\eps}
            h(x/\eps,u_\eps)\varphi\psi(x/\eps)
                        \dd\sigma_x\dd t\right|\\
 &\quad\leq\eps
   \left(\eps\int_0^T\int_{S_\eps}|h(x/\eps,u_\eps)|^2
                                       \dd\sigma_x\dd t\right)^{1/2}
   \left(\eps\int_0^T\int_{S_\eps}
                 |\varphi\psi(x/\eps)|^2\dd\sigma_x\dd t\right)^{1/2}
 \leq C\eps,
\end{align*}
by $(H_4)$, \eqref{eq:uniform-mass} and
\eqref{eq:surface-mass}. Letting $\eps\to0$ in
\eqref{eq:hom-oscillating-test} and using the two-scale flux
convergence \eqref{eq:twoscale-flux}, we obtain, for all such $\varphi$
and $\psi$,
\begin{equation}
 \int_Q\varphi\int_{\Ystar}A
       (D_xu_0+D_yu_1)\cdot D_y\psi\dd y\dd x\dd t=0.
 \label{eq:hom-cell-identified}
\end{equation}
For each fixed periodic smooth $\psi$, the inner integral in
\eqref{eq:hom-cell-identified} belongs to $L^1(Q)$, hence it
vanishes for a.e. $(x,t)$. Choosing a countable dense family of such
$\psi$ and using continuity with respect to the
$\mathcal V_\omega$ norm extends this identity to all
$z\in\mathcal V_\omega$. Subtracting the mean of $\psi$
does not change its gradient. Therefore $u_1(x,t,\cdot)$ solves
\eqref{eq:hom-cell} with $\xi=D_xu_0(x,t)$.
Proposition~\ref{prop:hom-cell} gives uniqueness and linearity,
proving \eqref{eq:hom-corrector-identification}.
Finally, integrate the two-scale limit
\eqref{eq:twoscale-flux} over $Y$ and substitute the identified
corrector. The definition \eqref{eq:hom-matrix} yields
\eqref{eq:hom-flux-convergence}.
\end{proof}

\subsection{Passage to the limit and the macroscopic equation}

With the corrector identified, the limit of the microscopic balance can be
computed, and uniqueness of the limit removes the subsequences. Recall
$\overline b$ and $\theta$ from $(H_3)$ and $h^{\rm hom}$ from
\eqref{eq:effective-reaction}; by $(H_4)$, $h^{\rm hom}$ is globally
Lipschitz and nondecreasing, with $h^{\rm hom}(0)=0$.

\begin{theorem}
\label{thm:hom-main}
Under $(H_1)$--$(H_5)$, let $u_\eps$ be the unique microscopic
solutions from Theorem~\ref{thm:micro-solvability}. There exists a
unique
\[
 u_0\in L^2(0,T;H_0^1(\Om))\cap H^1(0,T;H^{-1}(\Om))
\]
solving
\begin{equation}
 \left\{\begin{aligned}
  \overline b\,\partial_tu_0
      -\Div_x(A^{\rm hom}D_xu_0)+h^{\rm hom}(u_0)
       &=\theta f &&\text{in }Q,\\
  u_0&=0&&\text{on }\partial\Om\times(0,T),\\
  u_0(0)&=u^{\rm in}&&\text{in }L^2(\Om).
 \end{aligned}\right.
 \label{eq:hom-effective-PDE}
\end{equation}
The equation holds in $L^2(0,T;H^{-1}(\Om))$; equivalently, for a.e.
$t\in(0,T)$,
\begin{equation}
 \overline b\,\langle\partial_tu_0(t),v\rangle
 +\int_\Om A^{\rm hom}D_xu_0(t)\cdot Dv\dd x
 +\int_\Om h^{\rm hom}(u_0(t))v\dd x
 =\theta\int_\Om f(t)v\dd x
 \qquad\forall\,v\in H_0^1(\Om).
 \label{eq:hom-effective-weak}
\end{equation}
Moreover, the convergences
\eqref{eq:compact-strong}, \eqref{eq:compact-mass-strong},
\eqref{eq:twoscale-gradient-weighted},
\eqref{eq:twoscale-flux}, \eqref{eq:hom-flux-convergence},
\eqref{eq:reaction-limit} and \eqref{eq:dissipation-limit}, with $u=u_0$,
hold for the \emph{whole}
family as $\eps\to0$, with
\begin{equation}
 u_1(x,t,y)=\sum_{j=1}^N\pi^j(y)D_{x_j}u_0(x,t)
 \qquad\text{for a.e. }(x,t,y)\in Q\times\Ystar.
 \label{eq:hom-whole-corrector}
\end{equation}
\end{theorem}

\begin{proof}
Let $\phi\in C_c^\infty(\Om\times[0,T))$. Inserting $v=\phi|_{\Omeps}$ into
the microscopic weak formulation and integrating by parts in time gives
\begin{multline}
 -\int_Q m_\eps\partial_t\phi\dd x\dd t
 +\int_Q\chi_{\Omeps}A^\eps Du_\eps\cdot D_x\phi\dd x\dd t
 +\eps\int_0^T\int_{S_\eps}h(x/\eps,u_\eps)\phi\dd\sigma_x\dd t\\
   =\int_Q\chi_{\Omeps}f\phi\dd x\dd t
   +\int_\Om\rho_\eps u^{\rm in}\phi(0)\dd x;
 \label{eq:hom-micro-integrated}
\end{multline}
the last term is the initial mass $m_\eps(0)=\rho_\eps u^{\rm in}$, and
$\phi(T)=0$. By
\eqref{eq:periodic-material-average}, it converges to
$\overline b\int_\Om u^{\rm in}\phi(0)\dd x$.

The strong mass convergence \eqref{eq:compact-mass-strong}
implies
\[
 -\int_Qm_\eps\partial_t\phi\dd x\dd t
 \to-\overline b\int_Qu_0\partial_t\phi\dd x\dd t.
\]
Proposition~\ref{prop:hom-corrector} gives
\[
 \int_Q\chi_{\Omeps}A^\eps Du_\eps\cdot D_x\phi\dd x\dd t
 \to\int_Q A^{\rm hom}D_xu_0\cdot D_x\phi\dd x\dd t.
\]
Proposition~\ref{prop:surface-limit} applies with $u=u_0$, by
\eqref{eq:compact-material}, \eqref{eq:twoscale-u0-H1} and
\eqref{eq:micro-bound}; hence
\[
 \eps\int_0^T\int_{S_\eps}h(x/\eps,u_\eps)\phi
                                 \dd\sigma_x\dd t
 \to\int_Qh^{\rm hom}(u_0)\phi\dd x\dd t.
\]
Finally, \eqref{eq:periodic-material-average} with $a=1$ gives
\[
 \chi_{\Omeps}\rightharpoonup^*\theta
                  \quad\text{in }L^\infty(\Om),
 \qquad
 \int_Q\chi_{\Omeps}f\phi\dd x\dd t
        \to\theta\int_Qf\phi\dd x\dd t.
\]
Passing to the limit in \eqref{eq:hom-micro-integrated}, we obtain
\begin{multline*}
 -\overline b\int_Qu_0\partial_t\phi\dd x\dd t
 +\int_QA^{\rm hom}D_xu_0\cdot D_x\phi\dd x\dd t
 +\int_Qh^{\rm hom}(u_0)\phi\dd x\dd t\\
 =\theta\int_Qf\phi\dd x\dd t
 +\overline b\int_\Om u^{\rm in}\phi(0)\dd x.
\end{multline*}
Theorem~\ref{thm:twoscale-weighted} gives
$u_0\in L^2(0,T;H_0^1(\Om))$; Theorem~\ref{thm:strong-compactness}
gives $\partial_tu_0\in L^2(0,T;H^{-1}(\Om))$ and
$u_0(0)=u^{\rm in}$ in $H^{-1}(\Om)$. The Hilbert-triple continuity
$L^2(0,T;H_0^1(\Om))\cap H^1(0,T;H^{-1}(\Om))
\hookrightarrow C([0,T];L^2(\Om))$ promotes this to the asserted initial
condition in $L^2(\Om)$. Since $\phi$ is arbitrary, density yields
\eqref{eq:hom-effective-weak} and \eqref{eq:hom-effective-PDE}.

For uniqueness, let $u$ and $v$ be two solutions in the stated class,
set $z=u-v$, and subtract their weak equations. Testing by $z$, by means of the Hilbert-triple energy identity, gives
\[
 \frac{\overline b}{2}\|z(t)\|_{L^2(\Om)}^2
 +\int_0^t\int_\Om A^{\rm hom}Dz\cdot Dz\dd x\dd r
 +\int_0^t\int_\Om
       \bigl(h^{\rm hom}(u)-h^{\rm hom}(v)\bigr)(u-v)\dd x\dd r=0.
\]
All terms are nonnegative by $(H_3)$,
\eqref{eq:hom-ellipticity} and monotonicity of $h^{\rm hom}$.
Thus $z=0$, and the macroscopic solution is unique.
Existence follows from the subsequential limit just constructed.

It remains to remove subsequences. Every sequence $\eps_n\to0$
has a subsequence to which Theorems~\ref{thm:strong-compactness}
and \ref{thm:twoscale-weighted} apply. Its macroscopic limit must be
the unique solution of \eqref{eq:hom-effective-PDE}, and its
corrector must be \eqref{eq:hom-whole-corrector}, by
Proposition~\ref{prop:hom-corrector}. Therefore the usual
subsequence-contradiction argument gives strong convergence of the
whole family of extended states and masses and the asserted weak
two-scale convergences and weak convergence of physical fluxes.
The surface conclusions likewise hold for the whole family.
\end{proof}

\section{Energy convergence and a strong weighted corrector}
\label{sec:strong-corrector}

We strengthen the weak two-scale limit by comparing the microscopic
and macroscopic energies, first at each fixed time and then uniformly
on $[0,T]$. Since only the mass has a controlled time derivative,
the terminal storage term is treated by duality. The uniform-in-time
state convergence is then obtained in the storage-weighted norm.

\begin{lemma}
\label{lem:terminal-storage}
Under $(H_1)$--$(H_5)$, for every $t\in[0,T]$,
\begin{equation}
 \frac{\overline b}{2}\|u_0(t)\|_{L^2(\Om)}^2
 \leq \liminf_{\eps\to0}
 \frac12\|J_\eps u_\eps(t)\|_{H_{b,\eps}}^2.
 \label{eq:corrector-terminal-liminf}
\end{equation}
Here and below $u_0$ is the unique limit from
Theorem~\ref{thm:hom-main}; the microscopic energy is evaluated using its
continuous mass-space representative.
\end{lemma}

\begin{proof}
Fix $t\in[0,T]$ and $\zeta\in C_c^\infty(\Om)$. The elementary inequality
$|a|^2\geq 2az-|z|^2$, integrated with respect to
$b^\eps(x)\,\dd x$ on $\Omeps$, gives
\begin{equation}
 \frac12\|J_\eps u_\eps(t)\|_{H_{b,\eps}}^2
 \geq \langle m_\eps(t),\zeta\rangle
      -\frac12\int_\Om\rho_\eps(x)|\zeta(x)|^2\dd x.
 \label{eq:corrector-terminal-dual}
\end{equation}
On the right-hand side, $m_\eps(t)$ is the image of $J_\eps u_\eps(t)$
under $J_\eps^*$ followed by restriction to $H_0^1(\Om)$.
By \eqref{eq:compact-mass-strong},
$m_\eps(t)\to\overline b\,u_0(t)$ in $H^{-1}(\Om)$ for the whole family,
while \eqref{eq:compact-rho-average} yields
$\int_\Om\rho_\eps\zeta^2\to\overline b\int_\Om\zeta^2$.
Taking the lower limit in \eqref{eq:corrector-terminal-dual} gives
\[
 \liminf_{\eps\to0}\frac12\|J_\eps u_\eps(t)\|_{H_{b,\eps}}^2
 \geq \overline b\int_\Om u_0(t)\zeta\dd x
       -\frac{\overline b}{2}\int_\Om|\zeta|^2\dd x.
\]
Since $u_0\in C([0,T];L^2(\Om))$, taking the supremum over
$\zeta\in C_c^\infty(\Om)$ proves \eqref{eq:corrector-terminal-liminf}.
\end{proof}

For the corrector formula, set
\begin{equation*}
 \mathscr C_{ij}(y)=
 \begin{cases}
  \sqrt{w_i(y)}\bigl(\delta_{ij}+D_{y_i}\pi^j(y)\bigr),
        &y\in\Ystar,\\
  0,&y\in \mathcal O,
 \end{cases}
 \qquad 1\leq i,j\leq N,
\end{equation*}
and extend $\mathscr C$ $Y$-periodically. By
\eqref{eq:hom-cell-bound} and $(H_2)$,
$\mathscr C\in L^2_\#(Y)^{N\times N}$. Thus
$\mathscr C(y)D_xu_0(x,t)\in L^2(Q\times Y)^N$, but the expression
$\mathscr C(x/\eps)D_xu_0(x,t)$ need not be uniformly bounded
in $L^2(Q)^N$ for an arbitrary gradient $D_xu_0\in L^2(Q)^N$. We therefore approximate
the macroscopic gradient before evaluating the oscillatory corrector.

\begin{theorem}
\label{thm:strong-corrector}
Assume $(H_1)$--$(H_5)$. For $t\in[0,T]$, set
\[
 \mathfrak D_\eps(t)=\sum_{i=1}^N\int_0^t\int_{\Omeps}
       w_i^\eps|D_iu_\eps|^2\dd x\dd r,\qquad
 \mathfrak D_0(t)=\int_0^t\int_\Om A^{\rm hom}D_xu_0\cdot D_xu_0\dd x\dd r,
\]
$\mathfrak E_\eps(t)=\|J_\eps u_\eps(t)\|_{H_{b,\eps}}^2$ and
$\mathfrak E_0(t)=\overline b\|u_0(t)\|_{L^2(\Om)}^2$. Then, for the whole family,
\begin{equation}
 \lim_{\eps\to0}\sup_{0\leq t\leq T}
  \bigl(|\mathfrak D_\eps(t)-\mathfrak D_0(t)|+|\mathfrak E_\eps(t)-\mathfrak E_0(t)|\bigr)=0,
 \label{eq:corrector-uniform-energy}
\end{equation}
and the cumulative surface dissipation
$\eps\int_0^t\int_{S_\eps}h(x/\eps,u_\eps)u_\eps\dd\sigma_x\dd r$ converges to
$\int_0^t\int_\Om h^{\rm hom}(u_0)u_0\dd x\dd r$ uniformly on $[0,T]$.
Moreover, denoting by $U_\eps(t)$ the continuous mass-space
representative of $u_\eps$, viewed in $L^2(\Omeps,b^\eps\dd x)$, we have
\begin{equation}
 \lim_{\eps\to0}\sup_{0\leq t\leq T}
 \int_{\Omeps}b^\eps(x)|U_\eps(x,t)-u_0(x,t)|^2\dd x=0.
 \label{eq:corrector-uniform-weighted-state}
\end{equation}
Finally, with $g_{\eps,i}=\chi_{\Omeps}\sqrt{w_i^\eps}D_i u_\eps$ and
$G(x,t,y)=\mathscr C(y)D_xu_0(x,t)$, the weighted gradients $g_\eps$
converge strongly two-scale to $G$ in $L^2(Q)^N$. More concretely, if $v_\ell\in C_c^\infty(Q)$ satisfies
$v_\ell\to u_0$ in $L^2(0,T;H_0^1(\Om))$, then
\begin{equation}
 \lim_{\ell\to\infty}\lim_{\eps\to0}
 \left\|g_\eps-\mathscr C(x/\eps)D_xv_\ell\right\|_{L^2(Q)^N}=0.
 \label{eq:corrector-strong-physical}
\end{equation}
Consequently, there exists a choice $\ell=\ell(\eps)\to\infty$
for which the norm in \eqref{eq:corrector-strong-physical} tends to zero.
\end{theorem}

\begin{proof}
Fix $t\in(0,T]$ and set $Q_t=\Om\times(0,t)$ and
$Q_{\eps,t}=\Omeps\times(0,t)$. Strong bulk convergence and
\eqref{eq:periodic-material-average}, applied to $f u_0$ on $Q_t$,
give
\begin{equation}
 \int_{Q_{\eps,t}}fu_\eps\dd x\dd r
       \to\theta\int_{Q_t}fu_0\dd x\dd r.
 \label{eq:corrector-source}
\end{equation}
The conclusion \eqref{eq:dissipation-limit} of
Proposition~\ref{prop:surface-limit} also holds on $(0,t)$: its
hypotheses are inherited by restriction from $(0,T)$, and its proof
uses only estimates integrated over the chosen time interval. Hence
\begin{equation}
 \eps\int_0^t\int_{S_\eps}
       h(x/\eps,u_\eps)u_\eps\dd\sigma_x\dd r
  \to\int_{Q_t}h^{\rm hom}(u_0)u_0\dd x\dd r.
 \label{eq:corrector-surface-energy}
\end{equation}
Testing \eqref{eq:hom-effective-weak} with $u_0$ on $(0,t)$
(the Hilbert-triple identity applies on each such interval) gives
\begin{equation}
 \frac{\overline b}{2}\|u_0(t)\|_{L^2(\Om)}^2
  +\int_{Q_t}A^{\rm hom}D_xu_0\cdot D_xu_0
  +\int_{Q_t}h^{\rm hom}(u_0)u_0
  =\frac{\overline b}{2}\|u^{\rm in}\|_{L^2(\Om)}^2+\theta\int_{Q_t}fu_0.
 \label{eq:corrector-macro-energy}
\end{equation}
The initial energies also converge: since $|u^{\rm in}|^2\in L^1(\Om)$,
\eqref{eq:periodic-material-average} gives
\begin{equation}
 \frac12\int_{\Omeps}b^\eps|u^{\rm in}|^2\dd x
 =\frac12\int_\Om\rho_\eps|u^{\rm in}|^2\dd x
 \to\frac{\overline b}{2}\|u^{\rm in}\|_{L^2(\Om)}^2.
 \label{eq:corrector-initial-energy}
\end{equation}
Using the microscopic identity \eqref{eq:micro-energy} at this
same time $t$, the convergences
\eqref{eq:corrector-source}--\eqref{eq:corrector-initial-energy},
and the mass lower bound \eqref{eq:corrector-terminal-liminf},
we obtain
\begin{equation*}
 \limsup_{\eps\to0}\sum_{i=1}^N\int_{Q_{\eps,t}}
       w_i^\eps|D_i u_\eps|^2\dd x\dd r
   \leq\int_{Q_t}A^{\rm hom}D_xu_0\cdot D_xu_0\dd x\dd r.
\end{equation*}
For the reverse inequality, restrict the weak two-scale convergence
of $g_\eps$ to $Q_t$, that is, test it against functions supported
in $Q_t$. Weak lower semicontinuity on $Q_t\times Y$,
identification \eqref{eq:hom-corrector-identification}, and the cell
energy identity \eqref{eq:hom-energy-identity} imply
\begin{equation*}
 \int_{Q_t}A^{\rm hom}D_xu_0\cdot D_xu_0\dd x\dd r
 \leq\liminf_{\eps\to0}\sum_{i=1}^N\int_{Q_{\eps,t}}
             w_i^\eps|D_i u_\eps|^2\dd x\dd r.
\end{equation*}
The two bounds give $\mathfrak D_\eps(t)\to\mathfrak D_0(t)$ for $t>0$.
Comparing the microscopic and macroscopic energy identities once more, now
with the diffusion integrals known to converge, yields
$\mathfrak E_\eps(t)\to\mathfrak E_0(t)$ for the same $t$. Both assertions
also hold at $t=0$, since $\mathfrak D_\eps(0)=0$ and
$\mathfrak E_\eps(0)=\int_\Om\rho_\eps|u^{\rm in}|^2\dd x\to\mathfrak E_0(0)$ by
\eqref{eq:corrector-initial-energy}.

We now prove that the convergence is uniform in time. Besides the
diffusion energies $\mathfrak D_\eps$ and $\mathfrak D_0$ defined in the
statement, consider the cumulative surface dissipations
\[
 \mathfrak S_\eps(t)=\eps\int_0^t\int_{S_\eps}h(x/\eps,u_\eps)u_\eps
                         \dd\sigma_x\dd r,
 \qquad\mathfrak S_0(t)=\int_0^t\int_\Om h^{\rm hom}(u_0)u_0\dd x\dd r,
\]
and the source primitives
\[
 \mathfrak F_\eps(t)=\frac12\int_{\Omeps}b^\eps|u^{\rm in}|^2\dd x
     +\int_0^t\int_{\Omeps}fu_\eps\dd x\dd r,
 \qquad\mathfrak F_0(t)=\frac{\overline b}{2}\|u^{\rm in}\|_{L^2(\Om)}^2
     +\theta\int_0^t\int_\Om fu_0\dd x\dd r.
\]
Since $w_i\geq0$ and $h(y,r)r\geq0$, the functions $\mathfrak D_\eps$
and $\mathfrak S_\eps$ are nondecreasing. By the pointwise convergence
just proved and \eqref{eq:corrector-surface-energy}, they converge at every
$t\in[0,T]$ to $\mathfrak D_0$ and $\mathfrak S_0$, respectively. Both limit functions
are continuous, as their integrands belong to $L^1(Q)$. We follow
P\'olya's argument for monotone functions
\cite[Part~II, Problem~127]{PolyaSzego1972}. For any
$\kappa>0$, choose $0=t_0<\cdots<t_m=T$ with
$\mathfrak D_0(t_j)-\mathfrak D_0(t_{j-1})<\kappa$. Monotonicity gives, for
each $t\in[t_{j-1},t_j]$,
\[
 |\mathfrak D_\eps(t)-\mathfrak D_0(t)|
 \leq\max_{k\in\{j-1,j\}}|\mathfrak D_\eps(t_k)-\mathfrak D_0(t_k)|+\kappa.
\]
The maximum tends to zero because the partition is finite. Letting
$\kappa\to0$ proves $\|\mathfrak D_\eps-\mathfrak D_0\|_{C([0,T])}\to0$; the same
argument gives $\|\mathfrak S_\eps-\mathfrak S_0\|_{C([0,T])}\to0$.

The source primitives converge uniformly as well. Indeed,
\eqref{eq:corrector-source} and \eqref{eq:corrector-initial-energy} give
$\mathfrak F_\eps(t)\to \mathfrak F_0(t)$ for each $t$,
while Cauchy--Schwarz and \eqref{eq:micro-bound} yield
\[
 |\mathfrak F_\eps(t)-\mathfrak F_\eps(s)|
 \leq C\|f\|_{L^2(\Om\times(s,t))}\quad (0\leq s<t\leq T).
\]
The absolute continuity of $t\mapsto\int_0^t\int_\Om|f|^2\dd x\dd r$ makes this a
common modulus of continuity. A finite-partition argument therefore
gives $\|\mathfrak F_\eps-\mathfrak F_0\|_{C([0,T])}\to0$.
The microscopic and macroscopic energy identities now imply
\[
 \tfrac12\mathfrak E_\eps(t)
   =\mathfrak F_\eps(t)-\mathfrak D_\eps(t)-\mathfrak S_\eps(t),\qquad
 \tfrac12\mathfrak E_0(t)
   =\mathfrak F_0(t)-\mathfrak D_0(t)-\mathfrak S_0(t).
\]
Thus both energy convergences are uniform, proving
\eqref{eq:corrector-uniform-energy}. In particular, the surface
dissipation itself converges uniformly in cumulative time.

We prove the stronger state convergence
\eqref{eq:corrector-uniform-weighted-state}, using the identification
of $H_{b,\eps}$ with $L^2(\Omeps,b^\eps\dd x)$ from Section~\ref{sec:micro}.
By \eqref{eq:compact-rho-average} and continuity of
$u_0:[0,T]\to L^2(\Om)$,
\begin{equation}
 \sup_{0\leq t\leq T}\left|
 \int_\Om\rho_\eps|u_0(t)|^2\dd x
 -\overline b\|u_0(t)\|_{L^2(\Om)}^2\right|\to0.
 \label{eq:corrector-uniform-average}
\end{equation}
Indeed, the set $\{|u_0(t)|^2:0\leq t\leq T\}$ is compact
in $L^1(\Om)$; weak-star convergence of the uniformly bounded
$\rho_\eps$ is uniform on compact subsets of $L^1(\Om)$, by a
finite-net argument.

It remains to identify the mixed product uniformly. Let
$\Gamma_\mu=(\mathrm{Id}-\mu\Delta_D)^{-1}$ denote the Dirichlet resolvent
on $\Om$ and put $z_\mu(t)=\Gamma_\mu u_0(t)$. For each fixed
$\mu>0$, the function $z_\mu$ belongs to
$C([0,T];H_0^1(\Om))$, while
\[
 \lim_{\mu\to0}\sup_{0\leq t\leq T}
 \|z_\mu(t)-u_0(t)\|_{L^2(\Om)}=0.
\]
Here the latter uniform approximation follows from strong resolvent
convergence on $L^2(\Om)$ and compactness of
$\{u_0(t):t\in[0,T]\}$. The mass-energy bound
\eqref{eq:micro-bound}, boundedness of $b$, and
\eqref{eq:compact-mass-strong}, first with fixed $\mu$ and then
as $\mu\to0$, yield
\begin{equation}
 \sup_{0\leq t\leq T}\left|
 \int_{\Omeps}b^\eps U_\eps(t)u_0(t)\dd x
 -\overline b\|u_0(t)\|_{L^2(\Om)}^2\right|\to0.
 \label{eq:corrector-uniform-mixed}
\end{equation}
Indeed, the approximation error in the mixed integral is
bounded, uniformly in $\eps$ and $t$, by
$C\|u_0(t)-z_\mu(t)\|_{L^2(\Om)}$; for fixed $\mu$, the
remaining term is controlled by
\[
 \|m_\eps-\overline b u_0\|_{C([0,T];H^{-1}(\Om))}
 \sup_{0\leq t\leq T}\|z_\mu(t)\|_{H_0^1(\Om)}.
\]
The analogous limit-side approximation has the same bound.
Expanding the square in
\eqref{eq:corrector-uniform-weighted-state}, using the uniform
mass-energy convergence, \eqref{eq:corrector-uniform-average},
and \eqref{eq:corrector-uniform-mixed} proves the claim.

By Theorems~\ref{thm:twoscale-weighted} and \ref{thm:hom-main},
$g_\eps$ converges weakly two-scale to $G$ for the whole family.
The cell energy identity also gives
\[
 \int_Q\int_Y|G(x,t,y)|^2\dd y\dd x\dd t
    =\int_QA^{\rm hom}D_xu_0\cdot D_xu_0\dd x\dd t.
\]
Hence the convergence $\mathfrak D_\eps(T)\to\mathfrak D_0(T)$ reads
\begin{equation}
 \lim_{\eps\to0}\|g_\eps\|_{L^2(Q)^N}^2
       =\|G\|_{L^2(Q\times Y)^N}^2,
 \label{eq:corrector-strong-two-scale}
\end{equation}
which, together with the weak two-scale convergence, is the asserted strong
two-scale convergence.

Fix $\ell$. Since $D_xv_\ell$ is smooth and bounded and
$\mathscr C\in L^2_\#(Y)^{N\times N}$, the entries $\mathscr C_{ij}\mathscr C_{ik}$
belong to $L^1_\#(Y)$, and \eqref{eq:periodic-L1-average} gives
\begin{equation}
 \|\mathscr C(x/\eps)D_xv_\ell\|_{L^2(Q)}^2
 \to\int_Q\int_Y
       |\mathscr C(y)D_xv_\ell(x,t)|^2\dd y\dd x\dd t.
 \label{eq:corrector-test-norm}
\end{equation}
Since $D_xv_\ell$ is smooth and compactly supported,
\eqref{eq:periodic-multiplier} with $p=2$, applied componentwise to
$\mathscr C\in L^2_\#(Y)^{N\times N}$, gives
\begin{equation}
 \int_Qg_\eps\cdot\mathscr C(x/\eps)D_xv_\ell\dd x\dd t
 \to\int_Q\int_YG(x,t,y)\cdot
          \mathscr C(y)D_xv_\ell(x,t)\dd y\dd x\dd t.
 \label{eq:corrector-cross}
\end{equation}
Combining \eqref{eq:corrector-strong-two-scale},
\eqref{eq:corrector-test-norm} and \eqref{eq:corrector-cross},
we obtain the exact limit
\begin{equation}
 \lim_{\eps\to0}
 \|g_\eps-\mathscr C(x/\eps)D_xv_\ell\|_{L^2(Q)^N}^2
 =\int_Q\int_Y|\mathscr C(y)(D_xu_0-D_xv_\ell)|^2\dd y\dd x\dd t.
 \label{eq:corrector-exact-residual}
\end{equation}
The right-hand side is at most
$\|\mathscr C\|_{L^2(Y)^{N\times N}}^2\|D_xu_0-D_xv_\ell\|_{L^2(Q)^N}^2$,
which tends to zero as $\ell\to\infty$. This proves
\eqref{eq:corrector-strong-physical}. To obtain a single family,
choose increasing integers $\ell_k$ so that the right-hand side
of \eqref{eq:corrector-exact-residual} is at most $k^{-2}$.
For each $k$, fix $\eps_k$, decreasing to zero, such that the squared residual
on the left is at most $2k^{-2}$ whenever $0<\eps<\eps_k$.
Taking $\ell(\eps)=\ell_k$ for
$\eps_{k+1}\leq\eps<\eps_k$ completes the diagonal selection.
\end{proof}

The weighted convergence controls the states wherever the storage is
bounded below.

\begin{corollary}\label{cor:positive-storage}
Under $(H_1)$--$(H_5)$, for every $\lambda>0$,
\[
 \lim_{\eps\to0}\sup_{0\leq t\leq T}
 \int_{\Omeps\cap\{b^\eps\geq\lambda\}}|U_\eps(x,t)-u_0(x,t)|^2\dd x=0.
\]
\end{corollary}

\begin{proof}
On $\{b^\eps\geq\lambda\}$ the integrand is bounded by
$\lambda^{-1}b^\eps|U_\eps-u_0|^2$; apply
\eqref{eq:corrector-uniform-weighted-state}.
\end{proof}

\begingroup
\renewcommand{\bibliofont}{\fontsize{8}{9.5}\selectfont}
\providecommand{\bysame}{\leavevmode\hbox to3em{\hrulefill}\thinspace}
\providecommand{\MR}{\relax\ifhmode\unskip\space\fi MR }
\providecommand{\MRhref}[2]{%
  \href{http://www.ams.org/mathscinet-getitem?mr=#1}{#2}
}
\providecommand{\href}[2]{#2}

\endgroup

\end{document}